\documentclass{amsart}
\usepackage{amsmath,amssymb,amsthm,amscd,enumerate}
\usepackage[all]{xy}

\usepackage{color}

\def\Bbb{\mathbb}

\newcommand\bS{{}^b\kern-2ptS}

\newcommand\datver[1]{\def\datverp%
 {\par\boxed{\boxed{\text{Version: #1; Run: \today}}}}}
\datver{3.0D; Revised: June 21, 1999}

\newcommand\CC{\mathbb C}

\newcommand\RR{\mathbb R}
\newcommand\ZZ{\mathbb Z}

\newcommand\supp{\operatorname{supp}}

\newcommand\mP{\mathcal  P}

\newcommand\Aut{\operatorname{Aut}\nolimits}

\newcommand\wt[1]{\widetilde{#1}}
\newcommand\wh[1]{\widehat{#1}}
\newcommand\ov[1]{\overline{#1}}
\def\g{\gamma}
\def\E{{\mathcal E}}

\def\H{{\mathcal H}}

\def\ss{\subset}

\def\ti{{\operatorname{t\mbox{-}ind}\nolimits}}

\def\o{\otimes}

\def\Hom{\operatorname{Hom}\nolimits}

\def\a{\alpha}
\def\b{\beta}
\def\f{\varphi}
\def\s{\sigma}
\def\t{\tau}

\def\x{\chi}
\def\e{\varepsilon}
\def\La{\Lambda}
\def\la{\lambda}
\def\w{\wedge}

\def\F{\Phi}

\def\r{\rho}

\def\cV{{\mathcal V}}

\def\R{{\Bbb R}}
\def\C{{\Bbb C}}

\def\Lra{\Leftrightarrow}
\newcommand\ola[1]{\stackrel{#1}{\longrightarrow}}

\newcommand{\Mat}[4]{\left( \begin{array}{cc}
                            #1 & #2 \\
                            #3 & #4
                      \end{array} \right)}

\def\LS{longitudinally smooth }
\def\ve{{vert}}

\newtheorem{teo}{Theorem}[section]
\newtheorem{lem}[teo]{Lemma}
\newtheorem{prop}[teo]{Proposition}
\newtheorem{cor}[teo]{Corollary}

\theoremstyle{definition}
\newtheorem{dfn}[teo]{Definition}

\theoremstyle{remark}
\newtheorem{rmk}[teo]{Remark}

\newcommand\Id{\operatorname{Id}}
\newcommand\ie{{\em i.e.,} }

\newcommand\etimes{\boxtimes}

\newcommand\GR{\mathcal{G}}

\begin{document}
\title[Thom isomorphism]
{The Thom isomorphism in gauge-equivariant $K$-theory}

\author{Victor Nistor}
\address{Department of Mathematics,
Pennsylvania State University}
\email{nistor@math.psu.edu}
\urladdr{
http://www.math.psu.edu/nistor}
\author{Evgenij Troitsky}
\thanks{V. N.  was partially supported by NSF Grants DMS-9971951,
DMS 02-00808, and a ``collaborative research'' grant. E. T. was
partially supported by
RFFI Grant  05-01-00923, Grant for the support of
leading scientific schools 
and Grant ``Universities of Russia'' YP.04.02.530.
The present joint research
was started under the hospitality of MPIM (Bonn).}
\address{Dept. of Mech. and Math., Moscow State University,
119899 Moscow, Russia}
\email{troitsky@mech.math.msu.su}
\urladdr{
http://mech.math.msu.su/\~{}troitsky}


\begin{abstract}\
In a previous paper \cite{NT1}, we have introduced the
gauge-equivariant $K$-theory group $K_{\GR}^{0}(X)$ of a bundle
$\pi_{X} : X \to B$ endowed with a continuous action of a bundle
of compact Lie groups $p : \GR \to B$. These groups are the
natural range for the analytic index of a family of
gauge-invariant elliptic operators (\ie a family of elliptic
operators invariant with respect to the action of a bundle of
compact groups). In this paper, we continue our study of
gauge-equivariant $K$-theory. In particular, we introduce and
study products, which helps us establish the Thom isomorphism in
gauge-equivariant $K$-theory. Then we construct push-forward maps
and define the topological index of a gauge-invariant family.
\end{abstract}

\maketitle

\tableofcontents

\section{Introduction}

In this paper we establish a Thom isomorphism theorem for gauge
equivariant $K$-theory.  Let $p : \GR \to B$ be a bundle of {\em
compact} groups. Recall that this means that each fiber $\GR_b :=
p^{-1}(b)$ is a compact group and that, locally, $\GR$ is of the
form $U \times G$, where $U \subset B$ open and $G$ a fixed
compact group. Let $X$ and $B$ be locally compact spaces and
$\pi_{X} : X \to B$ be a continuous map. In the present
paper, as in \cite{NT1}, this map will be supposed to be a locally
trivial bundle. A part of present results can be
extended to the case of a general map, this,
as well as the proof of a general index theorem,
will be the subject of a forthcoming paper.

Assume that $\GR$ acts on $X$. This action will be always be
fiber-preserving. Then we can associate to the action of $\GR$ on
$X$ {\em $\GR$-equivariant $K$-theory groups} $K_{\GR}^{i}(X)$ as
in \cite{NT1}. We shall review and slightly generalize this
definition in Section \ref{sec.prel}.

For $X$ compact, the group $K_{\GR}^{0}(X)$ is defined as the
Grothendieck group of $\GR$-equivariant vector bundles on $X$. If
$X$ is not compact, we define the groups $K_{\GR}^{0}(X)$ using
fiberwise  one-point compactifications.
We shall call these groups simply {\em gauge-equivariant
$K$-theory groups} of $X$ when we do not want to specify $\GR$.
The reason for introducing the gauge-equivariant $K$-theory groups
is that they are the natural range for the index of a
gauge-invariant families of elliptic operators. In turn, the
motivation for studying gauge-invariant families and their index
is due to their connection to spectral theory and boundary value
problems on non-compact manifolds. Some possible connections with
Ramond-Ramond fields in String Theory were mentioned in
\cite{FreedCong,NT1}.
See also \cite{AStwist,FreHopTel,MaMeSi,NistorFam}.

In this paper, we continue our study of gauge-equivariant
$K$-theory. We begin by providing two alternative definitions of
the relative $K_{\GR}$--groups, both based on complexes of vector
bundles. (In this paper, all vector bundles are complex vector
bundles, with the exception of the tangent bundles and where
explicitly stated.) These alternative definitions, modeled on the
classical case \cite{AK, Friedrich}, provide a convenient
framework for the study of products, especially in the relative or
non-compact cases. The products are especially useful for the
proof of the Thom isomorphism in gauge-equivariant theory, which
is one of the main results of this paper. Let $E \to X$ be a
$\GR$-equivariant complex vector bundle. Then the Thom isomorphism
is a natural isomorphism
\begin{equation}\label{eq.intr.Thom}
    \tau_E : K^{i}_{\GR} (X) \to K^{i}_{\GR} (E).
\end{equation}
(There is also a variant of this result for $\rm spin^c$-vector
bundles, but since we will not need it for the index
theorem~\cite{NT3}, we will not discuss this in this paper.) The
Thom isomorphism allows us to define Gysin (or push-forward) maps
in $K$-theory. As it is well known from the classical work of
Atiyah and Singer \cite{AS1}, the Thom isomorphism and the Gysin
maps are some of the main ingredients used for the definition and
study of the topological index. In fact, we shall proceed along
the lines of that paper to define the topological index for
gauge-invariant families of elliptic operators. Some other approaches
to Thom isomorphism in general settings of Noncommutative
geometry were the subject of
\cite{ConnesThom,FaSk,HilsSkand,TroAcApMat2001,MaMeSi} and many other papers.

Gauge-equivariant $K$-theory behaves in many ways like the usual
equivariant $K$-theory, but exhibits also some new phenomena. For
example, the groups $K_{\GR}^{0}(B)$ may turn out to be reduced to
$K^{0}(B)$ when $\GR$ has ``a lot of twisting'' \cite[Proposition
3.6]{NT1}. This is never the case in equivariant $K$-theory when
the action of the group is trivial but the group itself is not
trivial. In \cite{NT1}, we addressed this problem in two ways:
first, we found conditions on the bundle of groups $p : \GR \to B$
that guarantee that $K^{0}_{\GR}(X)$ is not too small (this
condition is called {\em finite holonomy} and is recalled below),
and, second, we studied a substitute of $K^{0}_{\GR}(X)$ which is
never too small (this substitute is $K(C^*(\GR))$, the $K$-theory
of the $C^*$-algebra of the bundle of compact groups $\GR$).

In this paper, we shall again need the finite holonomy condition,
so let us review it now. To define the finite holonomy condition,
we introduced the {\em representation covering} of $\GR$, denoted
$\widehat \GR \to B$. As a space, $\widehat \GR$ is the union of
all the representation spaces $\widehat \GR_b$ of the fibers
$\GR_b$ of the bundle of compact groups $\GR$. One measure of the
twisting of the bundle $\GR$ is the holonomy associated to the
covering $\widehat \GR \to B$. We say that $\GR$ has {\em
representation theoretic finite holonomy} if $\widehat\GR$ is a
union of compact-open subsets. (An equivalent conditions can be
obtained in terms of the fundamental groups when $B$ is
path-connected, see Proposition~\ref{prop.r.f.h} below.)

Let $C^*(\GR)$ be the enveloping $C^*$-algebra of the bundle of
compact groups $\GR$. We proved in \cite[Theorem 5.2]{NT1} that
\begin{equation}\label{eq.isom.g.a}
        K^j_\GR(B) \cong K_j(C^*(\GR)),
\end{equation}
provided that $\GR$ has representation theoretic finite holonomy.
This guarantees that $K^j_\GR(B)$ is not too small. It also points
out an alternative, algebraic definition of the groups
$K_\GR^i(X)$.

The structure of the paper is as follows. We start from the
definition of gauge-equivariant $K$-theory and some basic results
from \cite{NT1}, most of them related to the ``finite holonomy
condition,'' a condition on bundles of compact groups that we
recall in Section \ref{sec.prel}. In Section \ref{sect:KiM} we
describe an equivalent definition of gauge-equivariant $K$-theory
in terms of complexes of vector bundles. This will turn out to be
especially useful when studying the topological index. In Section
\ref{subs:Thom} we establish the Thom isomorphism in
gauge-equivariant $K$-theory, and, in Section \ref{sec:gysinmaps},
we define and study the Gysin maps. The properties of the Gysin
maps allow us to define in Section \ref{sec:topindex} the
topological index and establish its main properties.

\section{Preliminaries\label{sec.prel}}

We now recall the definition of gauge-equivariant $K$-theory and
some basic results from \cite{NT1}.
An important part of our
discussion is occupied by the discussion of the finite holonomy condition
for a bundle of compact groups $p : \GR \to B$.

{\em All vector bundles in this paper are assumed to be {\bf
complex} vector bundles, unless otherwise mentioned and excluding
the tangent bundles to the various manifolds appearing below.}

\subsection{Bundles of compact groups and finite holonomy conditions}
We begin by introducing bundles of compact and locally compact
groups. Then we study finite holonomy conditions for bundles of
compact groups.

\begin{dfn}\label{def.blcg}\
Let $B$ be a locally compact space and let $G$ be a locally
compact group. We shall denote by $\Aut(G)$ the group of
automorphisms of $G$. A {\em bundle of locally compact groups
$\GR$ with typical fiber $G$ over $B$} is, by definition, a fiber
bundle $\GR \to B$ with typical fiber $G$ and structural group
$\Aut(G)$.
\end{dfn}

We fix the above notation. Namely, {\em from now on and throughout
this paper, unless explicitly otherwise mentioned, $B$ will be a
compact space and $\GR \to B$ will be a bundle of compact groups
with typical fiber $G$.}

We need now to introduce the representation theoretic holonomy of
a bundle of Lie group with compact fibers $p : \GR \to B$. Let
$\Aut (G)$ be the group of automorphisms of $G$.  By definition,
there exists then a principal $\Aut(G)$-bundle $\mP \to B$ such
that
\begin{equation*}
    \GR \cong \mP \times_{\Aut(G)} G : = (\mP \times G)/\Aut(G).
\end{equation*}
We shall fix $\mP$ in what follows.

Let $\widehat \GR$ be the (disjoint) union of the sets
$\widehat{\GR}_b$ of equivalence classes of irreducible
representations of the groups $\GR_b$. Using the natural action of
$\Aut(G)$ on $\wh G$, we can naturally identify $\widehat \GR$
with $\mP \times_{\Aut(G)} \wh G$ as fiber bundles over $B$.

Let $\Aut_0(G)$ be the connected component of the identity in
$\Aut (G)$. The group $\Aut_0(G)$ will act trivially on the set
$\wh G$, because the later is discrete. Let
\begin{equation*}
    H_R := \Aut(G) / \Aut_0(G), \quad \mP_0 := \mP/ \Aut_0(G),
    \;\; \text{ and } \;\; \widehat \GR \simeq
    \mP_0 \times_{H_R} \wh G.
\end{equation*}
Above, $\widehat \GR$ is defined because $\mP_0$ is an
$H_R$-principal bundle. The space $\widehat \GR$ will be called
the {\em representation space of $\GR$} and the covering $\widehat
\GR \to B$ will be called {\em the representation covering
associated to $\GR$}.

Assume now that $B$ is a path-connected, locally simply-connected
space and fix a point $b_0 \in B$. We shall denote, as usual, by
$\pi_1(B,b_0)$ the fundamental group of $B$. Then the bundle
$\mP_0$ is classified by a morphism
\begin{equation}
        \rho : \pi_1(B,b_0) \to H_R := \Aut(G)/ \Aut_0(G),
\end{equation}
which will be called {\em the holonomy of the representation covering of
$\GR$}.

For our further reasoning, we shall sometimes need the following
finite holonomy condition.

\begin{dfn}\label{def.c.f.h}\
We say that $\GR$ has {\em representation theoretic finite holonomy}
if every $\sigma \in \widehat{\GR}$ is contained in a compact-open
subset of $\widehat \GR$.
\end{dfn}

In the cases we are interested in, the above condition can be
reformulated as follows \cite{NT1}

\begin{prop}\label{prop.r.f.h}\
Assume that $B$ is path-connected and locally simply-connected.
Then $\GR$ has {\em representation theoretic finite holonomy} if,
and only if $\pi_1(B,b_0) \sigma \subset \widehat G$ is a finite
set for any irreducible representation $\sigma$ of $G$.
\end{prop}

{\it From now on we shall assume that $\GR$ has
representation theoretic finite holonomy.}

\subsection{Gauge-equivariant $K$-theory}
Let us now define the gauge equivariant $K$-theory groups of a
``$\GR$-fiber bundle'' $\pi_{Y} : Y \to B$. All our definitions
are well known if $B$ is reduced to a point (cf.
\cite{AK,Friedrich}). First we need to fix the notation.

If $f_i : Y_i \to B$, $i = 1, 2$, are two maps, we shall denote by
\begin{equation}
        Y_1 \times_B Y_2 := \{ (y_1,y_2) \in Y_1 \times Y_2,\,
        f_1(y_1) = f_2(y_2)\, \}
\end{equation}
their fibered product. Let $p : \GR \to B$ be a bundle of locally
compact groups and let $\pi_{Y} : Y \to B$ be a continuous map. We
shall say that $\GR$ {\em acts on $Y$} if each group $\GR_b$ acts
continuously on $Y_b := \pi^{-1}(b)$ and the induced map $\mu$
\begin{equation*}
        \GR \times_B Y :=\{ (g,y) \in \GR \times Y,\, p(g) =
        \pi_{Y}(y)\} \ni\, (g,y) \longrightarrow
        \mu(g,y) := gy \, \in Y
\end{equation*}
is continuous. If $\GR$ acts on $Y$, we shall say that $Y$
is {\em a $\GR$-space}. If, in addition to that, $Y \to B$ is also
locally trivial, we shall say that $Y$ is {\em a $\GR$-fiber
bundle}, or, simply, a {\em $\GR$-bundle}. This definition is a
particular case of the definition of the action of a
differentiable groupoid on a space.

Let $\pi_{Y} : Y \to B$ be a $\GR$-space, with $\GR$ a bundle of
compact groups over $B$. Recall that a vector bundle $\tilde
\pi_{E} : E \to Y$ is {\em a $\GR$-equivariant vector bundle} (or
simply {\em a $\GR$--equivariant vector bundle}) if
\begin{equation*}
    \pi_E := \pi_{Y} \circ \tilde \pi_{E} : E \to B
\end{equation*}
is a $\GR$-space, the projection
\begin{equation*}
    \tilde\pi_{E} : E_b := \pi_E^{-1}(b) \to Y_b := \pi_Y^{-1}(b)
\end{equation*}
is $\GR_b := p^{-1}(b)$ equivariant, and the induced action $E_y \to E_{gy}$
of $g \in \GR$, between the corresponding fibers of $E \to Y$, is linear for any
$y \in Y_b$, $g \in \GR_b$, and $b \in B$.

To define gauge-equivariant $K$--theory, we first recall some
preliminary definitions from \cite{NT1}. Let $\tilde\pi_E : E \to
Y$ be a $\GR$-equivariant vector bundle and let $\tilde\pi_{E'} :
E' \to Y'$ be a $\GR'$-equivariant vector bundle, for two bundles
of compact groups $\GR \to B$ and $\GR' \to B'$. We shall say that
$(\g,\f, \eta, \psi) : (\GR', E', Y', B') \to (\GR, E, Y, B)$ is a
 {\em $\gamma$--equivariant
morphism of vector bundles} if the following five conditions are
satisfied:
\begin{enumerate}[(i)]
\item\ $\g:\GR' \to\GR,$\  $\f:E' \to E,$\ $\eta : Y' \to Y,$\
and $\psi : B \to B',$
\item\ all the resulting diagrams are commutative,
\item\ $\f(ge)=\g(g) \f(e)$ for all $e\in E'_b$ and all $g \in \GR'_b,$,
\item\ $\g$ is a group morphism in each fiber, and
\item\ $f$ is a vector bundle morphism.
\end{enumerate}
We shall say that {\em $\phi : E \to E'$ is a
$\gamma$--equivariant morphism of vector bundles} if, by
definition, it is part of a morphism $(\g,\f, \eta, \psi) : (\GR',
E',Y', B') \to (\GR, E, Y, B)$. Note that $\eta$ and $\psi$ are
determined by $\gamma$ and $\phi$.

As usual, if $\psi : B' \to B$ is a continuous [respectively,
smooth] map, we define the {\em inverse image\/} $(\psi^*(\GR),
\psi^*(E), \psi^*(Y), B')$ of a $\GR$-equivariant vector bundle $E
\to Y$ by $\psi^*(\GR) = \GR \times_B B'$, $\psi^*(E) = E \times_B
B'$, and $\psi^*(Y) = Y \times_B B'$. If $B' \subset B$ and $\psi$
is the embedding, this construction gives the restriction of a
$\GR$-equivariant vector bundle $E \to Y$ to a closed, invariant
subset $B'\subset B$ of the base of $\GR$, yielding a
$\GR_{B'}$--equivariant vector bundle. Usually $\GR$ will be
fixed, however.

Let $p : \GR \to B$ be a bundle of compact groups and $\pi_Y : Y
\to B$ be a $\GR$--space. The set of isomorphism classes of
$\GR$--equivariant vector bundles $\tilde\pi_E : E \to Y$ will be
denoted by $\E_\GR(Y)$. On this set we introduce a monoid
operation, denoted ``$+$,'' using the direct sum of vector
bundles. This defines a monoid structure on the set $\E_\GR(Y)$ as
in the case when $B$ consists of a point.

\begin{dfn}\label{def.K}\
Let $\GR \to B$ be a bundle of compact groups acting on the
$\GR$-space $Y \to B$. Assume $Y$ to be compact. The {\em
$\GR$-equivariant $K$-theory group} $K^0_{\GR}(Y)$ is defined as
the group completion of the monoid $\E_\GR(Y)$.
\end{dfn}

When working with gauge-equivariant $K$--theory, we shall use the
following terminology and notation. If $E \to Y$ is a
$\GR$-equivariant vector bundle on $Y$, we shall denote by $[E]$
its class in $K^0_\GR(Y)$. Thus $K^0_\GR(Y)$ consists of
differences $[E] - [E^1]$. The groups $K^0_\GR(Y)$ will also be
called {\em gauge equivariant} $K$-theory groups, when we do not
need to specify $\GR$. If $B$ is reduced to a point, then $\GR$ is
group, and the groups $K^0_\GR(Y)$ reduce to the usual equivariant
$K$-groups.

We have the following simple observations on gauge-equivariant
$K$--theory. First, the familiar functoriality properties of the
usual equivariant $K$-theory groups extend to the gauge
equivariant $K$-theory groups. For example, assume that the bundle
of compact groups $\GR \to B$ acts on a fiber bundle $Y \to B$ and
that, similarly, $\GR' \to B'$ acts on a fiber bundle $Y' \to B'$.
Let $\g : \GR \to \GR'$ be a morphism of bundles of compact groups
and $f : Y \to Y'$ be a $\g$-equivariant map. Then we obtain a
natural group morphism
\begin{equation}
        (\g,f)^* : K^0_{\GR'}(Y') \to K^0_{\GR}(Y).
\end{equation}
If $\gamma$ is the identity morphism, we shall denote $(\g,f)^* =
f^*$.

A $\GR$-equivariant vector bundle $E \to Y$ on a $\GR$-space $Y
\to B$, $Y$ compact, is called {\em trivial} if, by definition,
there exists a $\GR$-equivariant vector bundle $E' \to B$ such
that $E$ is isomorphic to the pull-back of $E'$ to $Y$. Thus $E
\simeq Y \times_B E'$. If $\GR \to B$ has representation theoretic
finite holonomy and $Y$ is a compact $\GR$-bundle, then every
$\GR$--equivariant vector bundle over $Y$ can be embedded into a
trivial $\GR$--equivariant vector bundle. This embedding will
necessarily be as a direct summand.

If $\GR \to B$ does not have finite holonomy, it is possible to
provide examples if $\GR$--equivariant vector bundles that do not
embed into trivial $\GR$--equivariant vector bundles \cite{NT1}.
Also, a related example from \cite{NT1} shows that the groups
$K^0_{\GR}(Y)$ can be fairly small if the holonomy of $\GR$ is
``large.''

A further observation is that it follows from the definition that
the tensor product of vector bundles defines a natural ring
structure on $K^0_{\GR}(Y)$.  We shall denote the product of two
elements $a$ and $b$ in this ring by $a \otimes b$ or, simply,
$ab$, when there is no danger of confusion.  In particular the
groups $K^i_\GR(X)$ for $\pi_X : X \to B$ are equipped with a
natural structure of $K^0_\GR(B)$--module obtained using the
pull-back of vector bundles on $B$,
namely, $a b := \pi_X^*(a) \otimes b \in K^0_\GR (X)$ for $a \in
K^0_\GR(B)$ and $b \in K^0_\GR (X)$.

The definition of the gauge-equivariant groups extends to
non-compact $\GR$-spaces $Y$ as in the case of equivariant
$K$--theory. Let $Y$ be a $\GR$-bundle.  We shall denote then by
$Y^+ := Y \cup B$ the compact space obtained from $Y$ by the
one-point compactification of each fiber (recall that $B$ is
compact). The need to consider the space $Y^+$ is the main reason
for considering also non longitudinally smooth fibers bundles on
$B$. Then
\begin{equation*}
        K_{\GR}^0(Y) := \ker
        \big(K^0_\GR(Y^+) \to K^0_\GR(B)\big).
\end{equation*}
Also as in the classical case, we let
\begin{equation*}
    K^n_{\GR}(Y,Y') := K^0_{\GR}((Y \setminus Y') \times \R^n)
\end{equation*}
for a $\GR$-subbundle $Y'\ss Y$. Then \cite{NT1} we have the
following periodicity result

\begin{teo}\label{bott}\
We have natural isomorphisms
\begin{equation*}
     K^n_\GR(Y,Y') \cong K^{n-2}_\GR(Y,Y').
\end{equation*}
\end{teo}

Gauge-equivariant $K$-theory is functorial with respect to open
embeddings. Indeed, let $U \subset X$ be an open,
$\GR$-equivariant subbundle. Then the results of \cite[Section
3]{NT1} provide us with a natural map morphism
\begin{equation} \label{eq:conext}
        i_*:K^n_\GR(U) \to K^n_\GR(X).
\end{equation}
In fact, $i_*$ is nothing but the composition $K^n_\GR(U) \cong
K^n_\GR(X, X \setminus U) \to K^n_\GR(X).$

\subsection{Additional results} We now prove some more results on
gauge-equivariant $K$-theory.

Let $\GR \to B$ and $\H \to B$ be two bundles of compact groups
over $B$. Recall that an $\H$-bundle $\pi_X : X \to B$ is called
{\em free} if the action of each group $\H_b$ on the fiber $X_b$
is free (\ie $hx = x$, $x \in X_b$, implies that $h$ is the
identity of $\GR_b$.) We shall need the following result, which is
an extension of a result in \cite[page 69]{Friedrich}.
For simplicity, we shall write $\GR \times \H$ instead of $\GR \times_B \H$.

\begin{teo}\label{lem:superinduction}\
Suppose  $\pi_X : X \to B$ is  a $\GR\times\H$-bundle that is free
as an $\H$-bundle. Let $\pi : X \to X/\H$ be the (fiberwise)
quotient map. For any $\GR$--equivariant vector bundle
$\tilde{\pi}_E : E \to X/\H$, we define the {\em induced vector
bundle}
\begin{equation*}
        \pi^*(E) := \big \{(x,\e) \in X \times E,\, \pi(x) = \tilde
    \pi_E(\e) \big \} \to X,
\end{equation*}
with the action of $\GR \times \H$ given by
$(g,h)\cdot (x,\e) := ((g,h)x, g\e).$ Then $\pi^*$ is gives rise to
a natural isomorphism $K^0_{\GR}(X/\H) \to K^0_{\GR\times\H}(X).$
\end{teo}

\begin{proof}\
Let $\pi^* : K^0_{\GR}(X/\H) \to K^0_{\GR\times\H}(X)$ be the
induction map, as above. We will construct a map $r :
K^0_{\GR\times\H} (X)\to K^0_\GR(X/\H)$ satisfying $\pi^* \circ r
= \Id$ and $r \circ \pi^* = \Id$. Let $\pi_F : F \to X$ be a $\GR
\times \H$-vector bundle. Since the action of $\H$ on $X$ is free,
the induced map $\overline{\pi}_F : F/\H \to X/\H$ of quotient
spaces is a (locally trivial) $\GR$-bundle. Clearly, this
construction is invariant under homotopy, and hence we can define
$r[F]:=[F/\H]$.

Let us check now that $r$ is indeed an inverse of $\pi^*$. Denote by
$F \ni f \to \H f \in F/\H$ the quotient map. Let $F \to X$ be a
$\GR \times \H$-vector bundle. To begin with, the total space of
$\pi^* \circ r (F)$ is
\begin{equation*}
    \big\{(x, \H f)\in X\times (F/\H),\, \pi^*(x) =
    \overline{\pi}_F( \H f) \big\},
\end{equation*}
by definition. Then the map $F \ni f \to (\pi_F(f), \H f)
\in \pi^* \circ r (F)$ is an isomorphism. Hence, $\pi^* r = \Id$.

Next, consider a $\GR$--vector bundle $\pi_E : E \to X/\H$. The
total space of $r \circ \pi^* (E)$ is then
\begin{equation*}
        \big \{ (\H y, \e) \in (X/ \H) \times E ,\,
        \H y =\pi_E(\e)) \big \},
\end{equation*}
because $\H$ acts only on the first component of $\pi^*(E)$. Then
\begin{equation*}
    E \to r \circ \pi^* (E),\qquad \e \mapsto (\pi(\e), \e)
\end{equation*}
is an isomorphism. Hence, $r \circ \pi^*=\Id$.
\end{proof}

See also \cite[Theorem 3.5]{NT1}.

\begin{cor}\label{cor.ind}\
Let $P$ be a $\GR \times \H$--bundle that is free as an
$\H$--bundle. Also, let $W$ be an $\H$--bundle, then there is a
natural isomorphism
\begin{equation*}
    K_{\GR\times\H}(P\times_B W)\cong K_\GR(P\times_{\H}W).
\end{equation*}
\end{cor}

\begin{proof}\ Take $X := P \times_B W$ in the previous theorem.
Then $X$ is a free $\H$--bundle, because $P$ is, and $X/\H =: P
\times_\H W$.
\end{proof}

In the following section, we shall also need the following
quotient construction associated to a trivialization of a
vector bundle over a subset. Namely, if $Y \subset X$ is a
$\GR$--invariant, closed subbundle, then we shall
denote by $X/_BY$ the fiberwise quotient space over $B$, that is the
quotient of $X$ with respect to the equivalence relation $\sim$,
$x \sim y$ if, and only if, $x, y \in Y_b$, for some $b \in B$.

If $E$ is a $\GR$--equivariant vector bundle over a $\GR$--bundle
$X$, together with a $\GR$--equivariant trivialization over a
$\GR$--subbundle $Y \subset X$, then we can generalize the quotient
(or collapsing) construction of \cite[\S 1.4]{AK} to obtain a
vector bundle over $X/_B Y$, where by $X/_BY$ we denote the
fiberwise quotient bundle over $B$, as above.

\begin{lem}\label{lem:compressionGR}\
Suppose that $X$ is a $\GR$--bundle and that $Y\ss X$ is a closed,
$\GR$--invariant subbundle. Let $E\to X$ be a $\GR$--vector bundle
and $\a:E|_Y\cong Y\times_B \cV$ be a $\GR$--equivariant
trivialization, where $\cV\to B$ is a $\GR$--equivariant vector
bundle. Then we can naturally associate to $(E, \a)$ a naturally
defined vector bundle $E/\a\to X/_B Y$ that depends only on
homotopy class of $\a$.
\end{lem}

\begin{proof}\
Let $p:Y\times_B \cV \to \cV$ be the natural projection.
Introduce the following equivalence relation on $E$:
 \begin{equation*}
    e\sim e' \Lra e,e'\in E|_Y \mbox{ and } p\a(e)=p\a(e').
 \end{equation*}
Let then $E/\a$ be equal to $E/\sim$. This is locally trivial
vector bundle over $X/_B Y$. Indeed, it is necessary to verify
this only in a neighborhood of $Y/_B Y\cong B$. Let $U$ be a
$\GR$ invariant open subset of $X$ such that $\a$ can be extended
to an isomorphism $\wt\a:E|_U\cong U\times_B \cV$. We obtain an
isomorphism
 \begin{equation*}
    \a':(E|_U)/\a \cong (U/_B Y)\times_B \cV,\qquad
    \a'(e)=\wt\a (e).
 \end{equation*}
Moreover, $(U/_B Y)\times_B \cV$ is a locally trivial
$\GR$--equivariant vector bundle.

Suppose that $\a_0$ and $\a_1$ are homotopic trivializations
of $E|_Y$, that is,
trivializations such that there exists a trivialization
$\b:E\times I|_{Y\times I}\cong {Y\times I}\times_B \cV$,
$\b|_{E\times \{0\}}=\a_0$ and $\b|_{E\times \{1\}}=\a_1$.
Let
\begin{equation*}
    f : X/_B Y \times I \to (X\times I)/_{(B \times I)}
    (Y\times I).
\end{equation*}
Then the bundle $f^*((E\times I)/\b)$ over $X/_B Y \times I$
satisfies $f^*((E\times I)/\b)|_{(X/_B Y)\times \{i\}}=E/\a_i$,
$i=0,1$. Hence, $E/\a_0\cong E/\a_1$.
\end{proof}

 \section{$K$-theory and complexes}\label{sect:KiM}

For the purpose of defining the Thom isomorphism, it is convenient
to work with an equivalent definition of gauge-equivariant
$K$-theory in terms of complexes of vector bundles. This will turn
out to be especially useful when studying the topological index.

The statements and proofs of this section, except maybe Lemma
\ref{lem:comptoch1}, follow the classical ones \cite{AK,
Friedrich}, so our presentation will be brief.

\subsection{The $L^n_{\GR}$--groups}
We begin by adapting some well known concepts and constructions to
our settings.

Let $X \to B$ be a locally compact, paracompact
$\GR$--bundle. A {\em finite complex of $\GR$--equivariant vector
bundles} over $X$ is a complex
\begin{equation*}
    (E^*,d) = \Bigl( \dots \ola{d_{i-1}} E^i \ola{d_{i}} E^{i+1}
    \ola{d_{i+1}} \dots \Bigr)
\end{equation*}
$i \in \ZZ,$ of $\GR$--equivariant vector bundles over $X$ with
only finitely many $E^{i}$'s different from zero. Explicitly,
$E^i$ are $\GR$--equivariant vector bundles, $d_i$'s are
$\GR$--equivariant morphisms, $d_{i + 1} d_i=0$ for every $i$, and
$E^i=0$ for $|i|$ large enough. We shall also use the notation
$(E^*, d) =  \big(E^0, \ldots, E^{n},\; d_i : E^i\vert_Y \to
E^{i+1}\vert_Y\big),$, if $E^{i} = 0$ for $i < 0$ and for $i > n$.

As usual, a {\em morphism of complexes}\ $f: (E^*, d) \to (F^*,
\delta) $ is a sequence of morphisms $f_i: E^i\to F^i$ such that
$f_{i+1} d_{i} = \delta_{i+1} f_i,$ for all $i$. These
constructions yield the category of finite complexes of
$\GR$--equivariant vector bundles. Isomorphism in this category
will be denoted by $(E^*, d) \cong (F^*, \delta).$

In what follows, we shall consider a pair $(X,Y)$ of $\GR$--bundles
with $X$ is a compact $\GR$--bundle, unless explicitly otherwise
mentioned.

\begin{dfn}\label{def.LG}\
Let $X$ be a compact $\GR$--bundle and $Y$ be a closed
$\GR$--invariant subbundle. Denote by $C_\GR^n (X, Y)$ the set of
(isomorphism classes of) sequences
\begin{equation*}
    (E^*, d) = \big(E^0, E^{1}, \ldots, E^{n},\;
    d_k : E^k\vert_Y \to E^{k+1}\vert_Y\big)
\end{equation*}
of $\GR$--equivariant vector bundles over $X$ such that
$(E^k\vert_Y, d)$ is exact if we let $E^j = 0$ for $j < 0$ or $j> n$.

We endow $C_\GR^n (X, Y)$ with the semigroup structure given by
the direct sums of complexes. An element in $C_\GR^n(X,Y)$ is
called {\em elementary} if it is isomorphic to a complex of the
form \begin{equation*}
    \ldots \to 0 \to E \ola{\Id} E \to 0 \to \dots,
\end{equation*}
Two complexes $(E^*, d), (F^*, \delta) \in C^n_\GR(X,Y)$ are
called {\em equivalent\/} if, and only if, there exist elementary
complexes $Q^1,\dots, Q^k, P^1,\dots, P^m \in C^n_\GR(X,Y)$ such
that
\begin{equation*}
    E \oplus Q_1 \oplus \dots \oplus Q_k \cong F \oplus P_1
    \oplus \dots \oplus P_m.
\end{equation*}
We write $E \simeq F$ in this case. The semigroup of equivalence
classes of sequences in $C^n_\GR(X, Y)$ will be denoted by
$L_\GR^n(X,Y)$.
\end{dfn}

We obtain, from definition, natural injective semigroup
homomorphisms
\begin{equation*}
    C_\GR^n(X,Y)\to C_\GR^{n+1}(X,Y)\quad \text{and} \quad
    C_\GR(X,Y):=\bigcup_n C_\GR^n(X,Y).
\end{equation*}
The equivalence relation $\sim$ commutes with embeddings, so the
above morphisms induce morphisms $L_\GR^n(X,Y)\to
L_\GR^{n+1}(X,Y).$ Let $L_\GR^\infty (X,Y) :=
\displaystyle{\lim_\to}\, L_\GR^n(X,Y)$.

\begin{lem}\label{lem:prodolhom}\
Let $E \to X$ and $F \to X$ be $\GR$--vector bundles. Let $\a:
E\vert_Y\to F\vert_Y$ and $\b: E\to F$ be surjective morphisms of
$\GR$--equivariant vector bundles. Also, assume that $\a$ and
$\b|_Y$ are homotopic in the set of surjective $\GR$--equivariant
vector bundle morphisms. Then there exists a surjective morphism
of $\GR$--equivariant vector bundles $\widetilde\a : E \to F$ such
that $\widetilde \a|_Y=\a.$ The same result remains true if we
replace ``surjective'' with ``injective'' or ``isomorphism''
everywhere.
\end{lem}

\begin{proof}\
Let $Z:=(Y\times [0,1])\cup (X\times \{0\})$ and $\pi:Z\to X$ be
the projection. Let $\pi^*(E) \to Z$ and $\pi^* (F) \to Z$ be the
pull-backs of $E$ and $F$. The homotopy in the statement of the
Lemma defines a surjective morphism $a:\pi^*(E) \to \pi^*(F)$ such
that $a|_{Y\times \{1\}}=\a$ and $a|_{X\times \{0\}}=\b$. By
\cite[Lemma 3.12]{NT1}, the morphism $a$ can be extended to a
surjective morphism over $(U\times [0,1])\cup (X\times \{0\})$,
where $U$ is an open $\GR$--neighborhood of $Y$. (In fact, in that
Lemma we considered only the case of an isomorphism, but the case
of a surjective morphism is proved in the same way.) Let $\f:X\to
[0,1]$ be a continuous function such that $\f(Y)=1$ and
$\f(X\setminus U)=0$. By averaging, we can assume $\f$ to be
$\GR$--equivariant. Then define $\wt \a(x) = a(x,\f(x))$, for all
$x\in X$.
\end{proof}

\begin{rmk}\label{rem:ka25}\
Suppose that $X$ is a compact $\GR$--space and
$Y=\emptyset.$ Then we have a natural isomorphism $\x_1: L_\GR^1
(X,\emptyset) \to K^0_\GR (X)$ taking the class of $ E^1,E^0$ to
the element $ [E^0] - [E^1].$
\end{rmk}

We shall need the following lemma.

\begin{lem}\label{lem:comptoch1}\
Let $p : \GR \to B$ be a bundle of compact groups and $\pi_X : X
\to B$ be a compact $\GR$--bundle. Assume that $\pi_X$ has a
cross-section, which we shall use to identify $B$ with a subset of
$X$. Then the sequence
\begin{equation*}
    0\to L_\GR^1(X,B) \to L_\GR^1(X) \to L_\GR^1(B)
\end{equation*}
is exact.
\end{lem}

\begin{proof}\
Suppose that $E = (E^1, E^0, \f)$ defines an element of
$L_\GR^1(X)$ such that its image in $L_\GR^1(B)$ is zero. Then the
definition of $E \sim 0$ in $L^1_\GR(B)$ shows that the
restrictions of $E^1$ and $E^0$ to $B$ are isomorphic over $B$.
Hence, the above sequence is exact at $L_\GR^1(X)$.

Suppose now that $(E^1,E^0,\f)$ represents a class in
$L_\GR^1(X,B)$ such that its image in $L_\GR^1(X)$ is zero. This
means (keeping in mind Remark \ref{rem:ka25}) that there exists a
$\GR$--equivariant vector bundle $\wt P$ and an isomorphism
$\wt\psi : E^1 \oplus \wt P \cong E^0 \oplus \wt P$. Let us define
$P := \wt P \oplus \pi_X^*(E^0|_B) \oplus \pi_X^*(\wt P|_B)$,
where $\pi_X : X \to B$ is the canonical projection, as in the
statement of the Lemma. Also, define $\psi=\wt \psi \oplus \Id :
E^1 \oplus P \to E^0\oplus P$, which is also an isomorphism.

We thus obtain that $T := \psi(\f\oplus \Id)^{-1}$  is an
automorphism of $(E^0\oplus P)|_B$. which has the form $\Mat \b 00
{\Id}$ with the respect to the decomposition
\begin{equation*}
    (E^0\oplus P)|_B = (E^0\oplus \wt P)|_B \oplus
    (E^0\oplus \wt P)|_B.
\end{equation*}
The automorphism $T := \psi(\f\oplus \Id)^{-1}$ is homotopic to the
automorphism $T_1$ defined by the matrix
\begin{equation*}
    \Mat {\Id} 00\b.
\end{equation*}
Since $T_1$ extends to an automorphism of $E^0\oplus P$ over $X$,
namely $\Mat {\Id} 00 {\pi_X^*(\b)}$, Lemma \ref{lem:prodolhom}
gives that the automorphism $\psi(\f\oplus \Id)^{-1}$ also can be
extended to $X$, such that over $B$ we have the following
commutative diagram:
\begin{equation*}
    \xymatrix{
    (E^1\oplus P)|_B \ar[r]^{\f\oplus \Id} \ar[d]^{\psi|_B}
    & (E^0\oplus P)|_B \ar[d]^{\a|_B=
    {\Mat \b 00 {\Id}}}\\
    (E^0\oplus P)|_B  \ar[r]^{\Id}&  (E^0\oplus P)|_B.}
\end{equation*}
Hence, $(E^1,E^0,\f) \oplus (P,P,\Id) \cong (E^0 \oplus P,E^0
\oplus P,\Id)$ and so is zero in $L_\GR^1(X,B)$.
\end{proof}

\subsection{Euler characteristics}
We now generalize the above construction to other groups $L^n_\GR$, thus
proving the existence and uniqueness of Euler characteristics.

\begin{dfn}\label{def.euler.char}\ Let $X$ be a compact
$\GR$--space and $Y \subset X$ be a $\GR$--invariant subset.
An {\em Euler characteristic\/} $\x_n$ is a natural transformation
of functors $\x_n : L_\GR^n(X,Y) \to K^0_\GR(X,Y)$, such that for
$Y = \emptyset$ it takes the form
\begin{equation*}
    \x_n(E)=\sum_{i=0}^n (-1)^i [E^i],
\end{equation*}
for any sequence $E = (E^*, d) \in L^n_{\GR}(X, Y)$.
\end{dfn}

\begin{lem}\label{lem:eulerchar1}\
There exists a unique natural transformation of functors $($\ie an
Euler characteristic$)$
\begin{equation*}
    \x_1: L_\GR^1 (X, Y) \to K^0_\GR (X, Y),
\end{equation*}
which, for $Y = \emptyset$, has the form indicated in {\rm
\ref{rem:ka25}}.
\end{lem}

\begin{proof}\
To prove the uniqueness, suppose that $\chi_1$ and $\chi'_1$ are
two Euler characteristics on $L_\GR^1$. Then $\chi'_1\chi_1^{-1}$
is a natural transformation of $K_\GR^0$ that is equal to the
identity on each $K^0_\GR(X)$. Let us consider the long exact
sequence a long exact sequence
\begin{multline}\label{exact}
        \dots \to K^{n-1}_\GR(Y,Y')\to K^{n-1}_\GR(Y)\to
        K^{n-1}_{\GR_1}(Y') \\ \to K^n_\GR(Y,Y') \to K^n_\GR(Y) \to
        K^n_{\GR}(Y') \to \dots
\end{multline}
associated to a pair $(Y, Y)$ of $\GR$--bundles (see \cite{NT1},
Equation (10) for a proof of the exactness of this sequence). The
map $K^0_\GR(X,B)\to K^0_\GR(X)$ from this exact sequence is
induced by $(X,\emptyset)\to (X,B)$, and hence, in particular, it
is natural. Assume that $\pi_X : X \to B$ has a cross-section.
Then the exact sequence \eqref{exact} for $(Y, Y') = (X, B)$
yields a natural exact sequence $0 \to K^0_\GR(X, B)\to
K^0_\GR(X)$. This in conjunction with Lemma \ref{lem:comptoch1}
shows that $\chi'_1\chi_1^{-1}$ is the identity on $K^0_\GR(X,
B)$. Recall now that we agreed to denote by $X/_BY$ the fiberwise
quotient space over $B$, that is the quotient of $X$ with respect
to the equivalence relation $\sim$, $x \sim y$ if, and only if,
$x, y \in Y_b$, for some $b \in B$. Finally, since the map $(X,
Y)\to (X/_B Y, B)$ induces an isomorphism of $K^0_\GR$-groups
\cite[Theorem 3.19]{NT1}, $\chi'_1\chi_1^{-1}$ is the identity on
$K^0_\GR(X, Y)$ for all pairs $(X, Y)$.

To prove the existence of the Euler characteristic $\chi_1$, let
$(E^1,E^0,\a) := (\a : E^1 \to E^0)$ represent an element of
$L_\GR^1(X,Y)$. Suppose that $X_0$ and $X_1$ are two copies of $X$
and $Z:=X_0\cup_Y X_1 \to B$ is the $\GR$--bundle obtained by
identifying the two copies of $Y\ss X_i$, $i=0,1$. The
identification of $E^1|_Y$ and $E^0|_Y$ with the help of $\a$
gives rise to an element  $[F^0] - [F^1]\in
K^0_\GR(Z)$ defined as follows. By adding some bundle to both
$E^i$'s, we can assume that $E^1$ is trivial (that is, it is
isomorphic to the pull--back of a vector bundle on $B$). Then
$E^1$ extends to a trivial $\GR$--vector bundle $\wt E^1\to Z$. We
define $F^0 := E^0 \cup_\a E^1$ and $F^1 := \wt E^1$ .

The exact sequence \eqref{exact} and the natural $\GR$-retractions
$\pi_i:Z\to X_i$,
give natural direct sum decompositions
\begin{equation}\label{eq:decompforglue}
    K^0_\GR(Z) = K^0_\GR(Z,X_i) \oplus K^0_\GR(X_i), \qquad i=0,1.
\end{equation}
The natural map $(X_0,Y)\to (Z,X_1)$ induces an isomorphism
\begin{equation*}
    k : K^0_\GR (Z,X_1)\to K^0_\GR (X_0,Y).
\end{equation*}
Let us define then $\chi_1(E^0,E^1,\a)$ to be equal to the image
under $k$ of the $K^0_\GR(Z,X_1)$-component of $[E^1,\a, E^0]$
(with the respect to \eqref{eq:decompforglue}). It follows from
its definition that this map is natural, respects direct sums, and
is independent with respect to the addition of elementary
elements. Our proof is completed by observing that
$\chi_1(E^1,E^0,\a) = [E^0] - [E^1]$ when $Y = \emptyset$.
\end{proof}

We shall also need the following continuity property of the
functor $L^1_{\GR}$. Recall that we have agreed to denote by
$X/_BY$ the fiberwise quotient bundle over $B$.

\begin{lem}\label{lemma.iso.c}\
The natural homomorphism
\begin{equation*}
    \Pi^* : L_\GR^1(X/_B Y,Y/_B Y) = L_\GR^1(X/_B Y,B)
    \to L_\GR^1(X,Y)
\end{equation*}
is an isomorphism for all pairs $(X,Y)$ of compact $\GR$-bundles.
\end{lem}

\begin{proof}\
Lemmata \ref{lem:eulerchar1} and \ref{lem:comptoch1} give the
following commutative diagram
 \begin{equation*}
    \xymatrix{
    L_\GR^1(X/_B Y,B) \ar[r]^{\quad \Pi^*}\ar[d]_\cong^{\chi_1}&
    L_\GR^1(X,Y)\ar[d]^{\chi_1}\\
    K^0_\GR(X/_B Y,B) \ar[r]^{\quad \Pi^*}_{\quad  \cong} & K^0_\GR(X,Y).
    }
 \end{equation*}
From this we obtain injectivity.

To prove surjectivity, suppose that $E^1$ and $E^0$ are
$\GR$--equivariant vector bundles over $X$ and $\a:E^1|_Y\to
E^0|_Y$ is an isomorphism of the restrictions. Let $P\to X$ be a
$\GR$-bundle such that there is an isomorphism $\b:E^1\oplus
P\cong F$, where $F$ is a trivial bundle (\ie isomorphic to a pull
back from $B$). Then $(E^1,E^0,\a)\sim (F,E^0\oplus P,\g)$, where
$\g=(\a\oplus \Id)\b^{-1}$. The last object is the image of
$(F,(E^0\oplus P)/\g, \g/\g)$ (see Lemma~\ref{lem:compressionGR}).
\end{proof}

We obtain the following corollaries.

\begin{cor}\label{cor.isom.chi}\
The Euler characteristic $\x_1 : L_\GR^1(X,Y) \to K^0_\GR(X,Y)$ is
an isomorphism and hence it is defines an equivalence of functors.
\end{cor}

\begin{proof}\ This follows from Lemmas \ref{lemma.iso.c}
and \ref{lem:eulerchar1}.
\end{proof}

\begin{lem}\label{lem:homotinL}
The class of $(E^1,E^0,\a)$ in $L_\GR^1(X,Y)$ depends only on the
homotopy class of the isomorphism $\a$.
\end{lem}

\begin{proof}
Let $Z = X \times [0,1]$, $W = Y \times [0,1]$. Denote by $p : Z
\to X$ the natural projection and assume that $\a_t$ is a
homotopy, where $\a_0 = \a$. Then $\a_t$ gives rise an isomorphism
$\b : p^*(E^1)|_W\cong p^*(E^0)|_W$, and hence to an element
$(p^*(E^1), p^*(E^0), \b)$ of $L_\GR^1(Z,W)$. If
\begin{equation*}
    i_t : (X,Y)\to (X\times \{t\},Y\times \{t\})\ss (Z,W),
\end{equation*}
are the standard inclusions, then $(E^1,E^0,\a_t)=i^*_t(p^*(E^1),
p^*(E^0), \b)$. Consider the commutative diagram
\begin{equation*}
\xymatrix{
    L_\GR^1(X,Y) \ar[d]^{\chi_1}&
    L_\GR^1(Z,W) \ar[d]^{\chi_1}\ar[l]_{i^*_0}\ar[r]^{i^*_1}&
    L_\GR^1(X,Y) \ar[d]^{\chi_1}\\
    K^0_\GR(X,Y) & K^0_\GR(Z,W) \ar[l]_{i^*_0}\ar[r]^{i^*_1}&
    K^0_\GR(X,Y).
}
\end{equation*}
The vertical morphisms and and the morphisms of the bottom line of
the above diagram are isomorphisms. Hence, the arrows of the top
line are isomorphisms too. The composition $i^*_0 (i^*_1)^{-1}$ is
identity for the bottom line, hence it is the identity for the top
line too.
\end{proof}

The following theorem reduces the study of the functors
$L^n_{\GR}$, $n > 1$, to the study of $L^1_{\GR}$.

\begin{teo}\label{theorem.nat.iso}\
The natural map $j_n:L_\GR^n (X, Y) \to L_\GR^{n + 1} (X, Y)$ is
an isomorphism.
\end{teo}

\begin{proof}\
Let $E = (E^{0}, E^1, \ldots, E^{n+1}; d_k)$, $d_k : E^k\vert_Y
\to E^{k+1}\vert_Y$ represent an element of the semigroup
$L_\GR^{n+1}(X,Y)$. To prove the surjectivity of $j_n$, let us
first notice that $E$ is equivalent to the complex
\begin{multline*}
    (E^0,\, \dots,\, E^{n-2},\, E^{n-1}\oplus E^{n+1},\,
    E^n \oplus E^{n+1},\, E^{n+1};\\ d_0,\, \ldots,\,
    d_{n-2}\oplus 0,\, d_{n-1} \oplus \Id,\, d_{n}\oplus 0).
\end{multline*}
The maps\ $d_{n}\, \oplus\, 0 : (E^{n} \oplus E^{n+1})\vert_Y \to
E^{n+1}\vert_Y$ and\ $0\, \oplus\, \Id : (E^{n} \oplus
E^{n+1})\vert_Y \to E^{n+1}\vert_Y$ are homotopic within the set
of surjective, $\GR$-equivariant vector bundle morphisms $(E^{n}
\oplus E^{n+1})\vert_Y \to E^{n+1}\vert_Y$. Hence, by
Lemma~\ref{lem:prodolhom}, $d_{n} \oplus 0$ can be extended to a
surjective morphism $b : E^{n} \oplus E^{n+1} \to E^{n+1}$ of
$\GR$--equivariant vector bundles (over the whole of $X$). So, the
bundle $E^n \oplus E^{n+1}$ is isomorphic to $\ker (b) \oplus
E^{n+1}$. Hence, the $E$ is equivalent to
\begin{equation*}
    (E^0,\, \dots,\, E^{n-2},\, E^{n-1}\oplus E^{n+1},\,
    \ker(b),\, 0; d_0,\, \ldots,\, d_{n-2}\oplus 0,\,
    d_{n-1},\, 0).
\end{equation*}
This proves the surjectivity of $j_n$.

To prove the injectivity of $j_n$, it is enough to define, for any
$n$, a left inverse $q_n : L_\GR^{n}(X,Y) \to L_\GR^1(X,Y)$ to
$s_n := j_{n-1} \circ \dots \circ j_1$.  Suppose that $(E^{*}; d)$
represents an element of semigroup $L_\GR^{n}(X,Y)$. Choose
$\GR$-invariant Hermitian metrics on $E^i$ and let
$d^*_i:E^{i+1}|_Y\to E^i|_Y$ be the adjoint of $d_i$. Let
\begin{equation*}
    F^0:=\bigoplus_i E^{2i},\quad F^1:=\bigoplus_i
    E^{2i+1},\quad b : F^0|_Y\to F^1|_Y,\quad b = \sum_i (d_{2i} +
    d^*_{2i+1}).
\end{equation*}
A standard verification shows that $b$ is an isomorphism. Since
all invariant metrics are homotopic to each other, Lemma
\ref{lem:homotinL} shows that $(E,d)\to (F,b)$ defines a morphism
$q_n:L_\GR^n(X, Y) \to L_\GR^1(X, Y)$. This is the desired left
inverse for $s_n$.
\end{proof}

Let us observe that the proof of the above theorem and Lemma
\ref{lem:homotinL} give the following corollary.

\begin{cor}\label{cor:homotinL}\
The class of $E = (E^i, d_i)$ in $L^n_{\GR}(X, Y)$ does not change
if we deform the differentials $d_i$ continuously.
\end{cor}

We are now ready to prove the following basic result.

\begin{teo}
For each $n$ there exists a unique Euler characteristic
\begin{equation*}
    \chi_n:L_\GR^n(X,Y)\cong K^0_\GR (X,Y).
\end{equation*}
In particular, $L_\GR^\infty (X,Y)\cong K^0_\GR (X,Y)$ and $L_\GR^n(X, Y)$
has a natural group structure for any closed, $\GR$--invariant subbundle
$Y \subset X$.
\end{teo}

\begin{proof}\
The statement is obtained from the lemmas we have proved above as
follows. First of all, Theorem \ref{theorem.nat.iso} allows us to
define
\begin{equation*}
\chi_n := \chi_1 \circ j_1^{-1} \circ \ldots j_{n-1}^{-1} :
L^n_{\GR}(X, Y) \to K^0_{\GR}(X, Y).
\end{equation*}
Lemma \ref{lemma.iso.c} shows that $\chi_n$ is an isomorphism. The
uniqueness of $\chi_n$ is proved in the same way as the uniqueness
of $\chi_1$ (Lemma \ref{lem:eulerchar1}).
\end{proof}

\subsection{Globally defined complexes}
The above theorem provides us with an alternative definition of
the groups $K^0_{\GR}(X, Y)$. We now derive yet another definition
of these groups that is closer to what is needed in applications
and is based on differentials defined on $X$, not just on $Y$.

Let $(E, d)$ be a complex of $\GR$--equivariant vector bundles
over a $\GR$--space $X$.
A point $x\in X$ will be called a {\em point of
acyclicity of $(E,d)$\/} if the restriction of $(E,d)$ to $x,$ i.e.,
the sequence of linear spaces
\begin{equation*}
    (E,d)_x = \Bigl( \dots \ola{(d_i)_x} E^i_x \ola{(d_{i+1})_x}
    E^{i+1}_x \ola{(d_{i+2})_x} \dots \Bigr),
\end{equation*}
is exact. The {\em support $\supp (E, d)$}\ of the finite complex $(E, d)$
is the complement in $X$ of the set of its points of acyclicity. This
definition and the following lemma hold also for $X$ non-compact.

\begin {lem}\label{lemma.cl}\
The support $\supp (E,d)$ is a closed $\GR$--invariant subspace of $X.$
\end {lem}

\begin{proof}\
The fact that $\supp (E,d)$ is closed is  classical (see \cite{AK,Friedrich}
for example). The invariance should be checked up over one fiber of $X$ at
$b\in B$. But this is once again a well known fact of equivariant $K$-theory
(see e.g. \cite{Friedrich}).
\end{proof}

\begin{lem}\label{lem:prodocompl}\
Let $E^n, \dots, E^0$ be $\GR$--equivariant vector bundles over
$X$ and $Y$ be a closed, $\GR$--invariant subbundle of $X$. Suppose
there are given morphisms $d_i:E^i\vert_Y\to E^{i-1}\vert_Y$ such
that $(E^i\vert_Y, d_i)$ is an exact complex. Then the morphisms
$d_i$ can be extended to morphisms defined over $X$ such that we
still have a complex of $\GR$--equivariant vector bundles.
\end{lem}

\begin{proof}\ We will show that we can extend each $d_i$ to a morphism
$r_i : E^i \to E^{i - 1}$ such that $r_{i-1}\circ r_i = 0$. Let us
find a $\GR$-invariant open neighborhood $U$ of $Y$ in $X$ such
that for any $i$ there exist an extension $s_i$ of $d_i$ to $U$
with $(E,s)$ still an exact sequence. The desired $r_i$ will be
then be defined as $r_i=\r s_i$, where $\r:X\to [0,1]$ is a
continuous function, $\r =  1$ on $Y$ and $\supp \r \ss U$.

Let us construct $U$ by induction over $i$. Assume that for the
closure $\ov U_i$ of some open $\GR$-neighborhood of $Y$ in $X$ we can
extend $d_j$ to $s_j$ ($j=1,\dots,i$) such that on $\ov U_i$ the
sequence
 \begin{equation*}
    E^i\ola{s_i}E^{i-1}\ola{s_{i-1}}\dots\to E^0\to 0
 \end{equation*}
is exact. Suppose, $K_i := \ker (s_i|{\ov U_i})$. Then $d_{i+1}$
determines a cross section of the bundle $\Hom_\GR (E^i,K_i)|_Y$.
This section can be extended to an open $\GR$-neighborhood $V$ of
$Y$ in $\ov U_i$. We hence obtain an extension $s_{i+1}:E^{i+1}\to
K_i$ of $d_{i+1}:E^{i+1}\to K_i$ over $V$. Since $d_{i+1}|_Y$ is
surjective (with range $K_i$), the morphism $s_{i+1}$ will be
surjective $\ov U_{i+1}$ for some open $U_{i+1}\ss U_{i}$.
\end{proof}

The above lemma suggests the following definition.

\begin{dfn}\label{def.new.L}\ Let $X$ be a compact $\GR$--bundle
and $Y \subset X$ be a $\GR$--invariant subbundle. We define
$E_\GR^n(X,Y)$ to be the semigroup of homotopy classes of
complexes of $\GR$--equivariant vector bundles of length $n$ over
$X$ such that their restrictions to $Y$ are acyclic (\ie exact).

We shall say that two complexes are {\em homotopic} if they are
isomorphic to the restrictions to $X \times \{0\}$ and $X \times
\{1\}$ of a complex defined over $X\times I$ and acyclic over
$Y\times I$.
\end{dfn}

\begin{rmk}\label{remark.new}\
By Corollary~\ref{cor:homotinL}, the restriction of morphisms
induces a morphism $\F^n : E_\GR^n(X,Y)\to L_\GR^n(X,Y)$.
\end{rmk}

\begin{teo}\label{teo:isomLD}\ Let $X$ be a compact $\GR$--bundle
and $Y \subset X$ be a $\GR$--invariant subbundle. Then the natural
transformation $\F_n$, defined in the above remark, is an isomorphism.
\end{teo}

\begin{proof}
The surjectivity of $\F_n$ follows from \ref{lem:prodocompl}. The
injectivity of $\F_n$ can be proved in the same way as~\cite[Lemma
2.6.13]{AK}, keeping in mind Lemma \ref{lem:prodocompl}.

More precisely, we need to demonstrate that differentials of any
complex over $\GR$-subbundle $(X\times \{0\})\cup (X\times \{1\})
\cup (Y\times I)$ of $X\times I$, which is acyclic over $Y\times
I$, can be extended to a complex over the entire $X\times I$. The
desired construction has the following three stages. First, let
$V$ be a $\GR$-invariant neighborhood of $Y$ such that the
restriction of our complex is still acyclic on $(V\times
\{0\})\cup (V\times \{1\}) \cup (Y\times I)$  as well as on its
closure $(\bar V\times \{0\})\cup (\bar V\times \{1\}) \cup
(Y\times I)$. By Lemma \ref{lem:prodocompl}, one can extend the
differentials $d_i$ to $\GR$--equivariant morphisms $r_i$ over
$\bar V \times I$ that still define a complex. Second, let $\r_1$,
$\r_2$ be a $\GR$-invariant partition of unity subordinated to
covering $V\times I$, $(X\setminus Y)\times I$ of $X\times Y$. Let
us extend original differentials to $(X\times [0, 1/4])\cup
(X\times [3/4, 1]) \cup (V\times I)$ by taking
$d_i(x,t):=d_i(x,0)$ for $t\le \frac 14\cdot \r_2(x)$, $x\in
X\setminus Y$. Similarly near $t=1$. Also,
\begin{equation*} d_i(x,t):=r_i\left(x,\frac{\left(t-\frac 14\cdot
\r_2(x)\right)}{\left(1-\frac 12\cdot \r_2(x)\right)}\right),\quad
x\in V. \end{equation*} Third, multiplying this differential by a
function $\tau:X\times I\to I$ equal to $1$ on the original subset
of definition of the differential and to $0$ off $(X\times [0,
1/4])\cup (X\times [3/4, 1]) \cup (V\times I)$ we obtain the
desired extension.
\end{proof}

\subsection{The non-compact case}
In the case of a locally compact, paracompact $\GR$--bundle $X$, we
change the definitions of $L_\GR^n$ and $E_\GR^n$ as follows. In
the definition of $L_\GR^n$, the morphisms $d_i$ have to be
defined and form an exact sequence off the interior of some
compact $\GR$--invariant subset $C$ of $X \smallsetminus Y$ (the
complement of $Y$ in $X$). In the definition of $E_\GR^n$, the
complexes have to be exact outside some compact $\GR$--invariant
subset of $X \smallsetminus Y$. In other words, $L_\GR^n(X, Y)
=L_{\GR}^n(X^{+}, Y^{+})$.

Since the proof of Lemma~\ref{lem:prodocompl} is still valid, we
have the analogue of Theorem~\ref{teo:isomLD}: there is a natural
isomorphism
\begin{equation*}
    L_\GR^n(X,Y)\cong E_\GR^n(X,Y).
\end{equation*}

The proof of the other statement also can be extended to the
non-compact case. The only difference is that we have to replace
$Y$ with $X \setminus U$, where $U$ is an open, $\GR$--invariant
subset with compact closure. Then, when we study two element
sequences $E = (E^i, d_i)$, we have to take the unions of the
corresponding open sets. Of course, these sets are not bundles,
unlike $Y$, but for our argument using extensions this is not
dangerous. This ultimately gives
\begin{equation}\label{ur:izomLKnonc}
    K^0_\GR (X,Y) \cong L_\GR^n(X,Y) \cong E_\GR^n(X,Y),
    \quad n \ge 1.
\end{equation}
As we shall see below, the liberty of using these equivalent
definitions of $K^0_\GR(X, Y)$ is quite convenient in
applications, especially when studying products.

\section{The Thom isomorphism}\label{subs:Thom}

In this section, we establish the Thom isomorphism in
gauge-equivariant $K$-theory. We begin with a discussion of
products and of the Thom morphism.

\subsection{Products}
Let $\pi_X : X \to B$ be a $\GR$--space, $\tilde\pi_F : F \to X$
be a complex $\GR$--vector bundle over $X$, and $s: X \to F$ a
$\GR$--invariant section. We shall denote by $\La^i F$ the $i$-th
exterior power of $F$, which is again a complex $\GR$--equivariant
vector bundle over $X$. As in the proof of the Thom isomorphism
for ordinary vector bundles, we define the complex $\La (F,s)$ of
$\GR$--equivariant vector bundles over $X$ by
\begin{equation}\label{eq.comp.s}
        \La(F,s):=(0\to \La^0F \ola{\a^0} \La^1F \ola{\a^1} \dots
        \ola{\a^{n-1}}\La^nF\to 0),
\end{equation}
where $\a^k(v_x) = s(x) \w v_x$ for $v_x\in\La^kF^x$ and $n = \dim
F$. It is immediate to check that $\alpha^{j+1}(x) \alpha^{j}(x) =
0$, and hence that $(\La (F,s), \a)$ is indeed a complex.

The K\"unneth formula shows that the complex $\La(F, s)$ is
acyclic for $s(x) \not = 0$, and hence $\supp (\La(F,s)) := \{x\in
X|s(x)=0\}$. If this set is compact, then the results of Section
\ref{sect:KiM} will associate to the complex $\La(F, s)$ of
Equation \eqref{eq.comp.s} an element
\begin{equation}\label{eq.elem.s}
        [\La(F,s)]\in K^0_\GR(X).
\end{equation}

Let $X$ be a $\GR$--bundle and $\pi_F : F \to X$ be a
$\GR$--equivariant vector bundle over $X$. The point of the above
construction is that $\pi_{F}^*(F)$, the lift of $F$ back to
itself, has a canonical section whose support is $X$.  Let us
recall how this is defined. Let $\pi_{FF}: \pi_{F}^* (F) \to F$ be
the $\GR$--vector bundle over $F$ with total space
 \begin{equation*}
        \pi_{F}^* (F) := \{(f_1,f_2) \in F \times F,\,
        \pi_{F}(f_1) = \pi_{F}(f_2)\}
 \end{equation*}
and $\pi_{FF} (f_1,f_2) = f_1$. The vector bundle $\pi_{FF} :
\pi_{F}^* (F) \to F$ has the canonical section
\begin{equation*}
        s_F : F \to \pi_{F}^{*}F,\qquad s_F(f) = (f,f).
\end{equation*}
The support of $s_F$ is equal to $X$. Hence, if $X$ is a compact
space, using again the results of Section \ref{sect:KiM},
especially \ref{teo:isomLD}, we obtain an element
\begin{equation}\label{eq.sect.Thom}
        \la_F := [\La(\pi_F^* (F) ,s_F)] \in K^0_\GR(F).
\end{equation}

Recall that the tensor product of vector bundles defines a natural
product $ab = a \otimes b\in K^0_{\GR}(X)$ for any $a \in
K^0_{\GR}(B)$ and any $b \in K^0_{\GR}(X)$, where $\pi_X : X \to
B$ is a compact $\GR$--space, as above.

Recall that all our vector bundles are assumed to be complex
vector bundles, except for the ones coming from geometry (tangent
bundles, their exterior powers) and where explicitly mentioned.
Due to the importance that $F$ be complex in the following
definition, we shall occasionally repeat this assumption.

\begin{dfn}\label{ind:Thomhomomorphism}\
Let $\pi_F : F \to X$ be a (complex) $\GR$--equivariant
vector bundle. Assume the $\GR$--bundle $X \to B$ is compact and
let $\la_F \in K^0_\GR(F)$ be the class defined in Equation
\eqref{eq.sect.Thom}, then the mapping
\begin{equation*}
    \f^F : K^0_\GR(X)\to K^0_\GR(F),\qquad \f^F(a)
    = \pi_F^*(a) \otimes \la_F.
\end{equation*}
is called the {\em Thom morphism.}
\end{dfn}

As we shall see below, the definition of the Thom homomorphism
extends to the case when $X$ is not compact, although the Thom
element itself is not defined if $X$ is not compact.

The definition of the Thom isomorphism immediately gives the
following proposition. We shall use the notation of Proposition
\ref{ind:Thomhomomorphism}.

\begin{prop}\label{prop:thom3}\ The Thom morphism
$\f^F : K^0_\GR(X)\to K^0_\GR(F)$ is a morphism of
$K_{\GR}^0(B)$--mo\-du\-les.
\end{prop}

Let $\iota : X \hookrightarrow F$ be the zero section embedding of
$X$ into $F$. Then $\iota$ induces homomorphisms
\begin{equation*}
    \iota^* : K^0_\GR(F) \to K^0_\GR(X) \quad \text{and}
    \quad \iota^* \circ \f^F : K^0_\GR(X) \to K^0_\GR(X).
\end{equation*}
It follows from the definition that $\iota^* \f^F(a) = a \cdot
\sum_{i=0}^n(-1)^i \La^i F.$
\medskip

\subsection{The non-compact case}
We now consider now the case when $X$ is locally compact, but not
necessarily compact. The complex $\La (\pi_F^* (F), s_F)$ has a
non-compact support, and hence it does not define an element of
$K^0_\GR (F)$. However, if $a = [(E, \a)] \in K^0_\GR (X)$ is
represented by the complex $(E, \a)$ of vector bundles with
compact support (Section \ref{sect:KiM}), then we can still
consider the tensor product complex
\begin{equation*}
    \big(\pi_F^*(\E),\,\pi_F^*(\a) \big) \o \La(\pi_{F}^{*}F,s_F).
\end{equation*}
From the K\"unneth formula for the homology of a tensor product we
obtain that the support of a tensor product complex is the
intersection of the supports of the two complexes. In particular,
we obtain
\begin{multline}\label{eq.Ncomp.p}
    \supp \{ (\pi_{F}^*E, \pi_F^*\a) \o
    \La(\pi_{F}^{*}F, s_F) \} \ss
    \supp (\pi_{F}^*E, \pi_F^*\a) \cap
    \supp \La(\pi_{F}^{*}F, s_F) \ss \\
    \ss \supp (\pi_{F}^*E, \pi_F^*\a)
    \cap X=\supp (E, \a).
\end{multline}
Thus, the complex $(\pi_F^*\E,\pi_E^*\a)\o \La(\pi_{F}^{*}F,s_F)$
has compact support and hence defines an element in
$K^0_{\GR}(F)$.

\begin{prop}\label{prop.Ncomp}\
The homomorphism of $K^0_\GR(B)$-modules
\begin{equation}\label{eq.Ncompact}
        \f^{F} : K^0_\GR(X) \to K^0_\GR(F), \quad \f^{F}(a) =
        [(\pi_{F}^*\E, \pi_F^*\a) \o \La(\pi_{F}^{*}F, s_F)],
\end{equation}
defined in Equation \eqref{eq.Ncomp.p} extends the Thom morphism
to the case of not necessarily compact $X$. The Thom morphism
$\f^{F}$ satisfies
\begin{equation}
        i^* \f^{F}(a) = a \cdot \sum_{i=0}^n(-1)^i \La^i F
\end{equation}
in the non-compact case as well.
\end{prop}

Let $F \to X$ be a $\GR$--equivariant vector bundle and $F^1 = F
\times \RR$, regarded as a vector bundle over $X \times \RR$. The
periodicity isomorphisms in gauge-equivariant $K$-theory groups
\cite[Theorem 3.18]{NT1}
\begin{equation*}
    K^{i \pm 1}_{\GR}(X \times \RR, Y \times \RR)
    \simeq  K^{i}_{\GR}(X, Y)
\end{equation*}
can be composed with $\f^{F^1}$, the Thom morphism for $F^1$,
giving a morphism
\begin{equation}
        \f^F : K^{i}_{\GR}(X) \to K^{i}_{\GR}(F)\,, \quad i=0,1.
\end{equation}
This morphism is the Thom morphism for $K^1$.

Let $p_X : X \to B$ and $p_Y : Y \to B$ be two {\em compact}
$\GR$-fiber bundles. Let $\pi_E : E \to X$ and $\pi_F : F \to Y$
be two complex $\GR$--equivariant vector bundles. Denote by $p_1 :
X \times_{B} Y \to X$ and by $p_2 : X \times_{B} Y \to Y$ the
projections onto the two factors and define $E \etimes F :=
p_{1}^{*}E \otimes p_{2}^{*}F$. The $\GR$--equivariant vector
bundle $E \etimes F$ will be called the {\em external tensor
product of $E$ and $F$ over $B$}. It is a vector bundle over $X
\times_B Y$. Then the formula
\begin{equation*}
        K_{\GR}^0(X) \otimes K_{\GR}^0(Y) \ni [E] \otimes [F]
        \to [E] \etimes [F] := [E \etimes F] \in K_{\GR}^{0}(X
        \times_{B} Y)
\end{equation*}
defines a product $K_{\GR}^0(X) \otimes K_{\GR}^0(Y) \to
K_{\GR}^{0}(X \times_{B} Y)$.

In particular, consider two complex $\GR$--equivariant vector
bundles $\pi_E : E \to X$ and $\pi_F : F \to X$. Then $E \oplus F
= E \times_X F$ and we obtain a product
 \begin{equation*}
        K_{\GR}^{i}(E) \o K_{\GR}^{j}(F) \ni [E] \otimes [F]
        \to [E] \etimes [F] := [E \etimes F] \to K_\GR^{i+j}(E
        \oplus F).
 \end{equation*}
Using also periodicity, we obtain the product
\begin{equation}\label{eq.ext.prod}
        \etimes : K^{i}_{\GR}(E) \o K^{j}_{\GR}(F) \to
        K^{i + j}_{\GR}(E \oplus F).
\end{equation}
This product is again seen to be given by the tensor product of
the (lifted) complexes (when representing $K$-theory classes by
complexes) as in the classical case.

The external product $\etimes$ behaves well with respect to the
``Thom construction,'' in the following sense. Let $F^1$ and $F^2$
be two complex bundles over $X,$ and $s_1$, $s_2$ two
corresponding sections of these bundles. Then
\begin{equation}\label{eq.mult.Thom0}
        \La(F^1\oplus F^2, s_1\otimes 1 + 1 \otimes s_2)
        = \La(F^1, s_1) \etimes \La(F^2, s_2).
\end{equation}
In particular, if $X$ is compact, we obtain
\begin{equation}\label{eq.mult.Thom}
    \la_E \etimes \la_F = \la_{E\oplus F}.
\end{equation}
We shall write $s_1 + s_2 = s_1\otimes 1 + 1 \otimes s_2$, for
simplicity.

The following theorem states that the Thom class is multiplicative
with respect to direct sums of vector bundles (see also
\cite{Carvalho}).

\begin{teo}\label{thm:thom6}\
Let $E, F\to X$ be two $\GR$--equivariant vector bundles, and
regard $E \oplus F \to E$ as the $\GR$--equivariant vector bundle
$\pi_{E}^*(F)$ over $E$. Then $\f^{\pi_{E}^*(F) } \circ \f^{E} =
\f^{E\oplus F}.$
\end{teo}

The above theorem amounts to the commutativity of the diagram
\begin{equation}
    \xymatrix{  K^*_\GR(X) \ar[rr]^{\f^{E}}\ar[dr]_{\f^{E\oplus F}}
    & & K^*_\GR(E) \ar[dl]^{\f^{\pi_{E}^*(F) }}\\ & K^*_\GR(E\oplus F)
    & }
\end{equation}

\begin{proof}\
Let $F^1 := \pi_E^*(F) = E \oplus F$, regarded as a vector bundle
over $E$. Consider the projections
\begin{equation*}
    \pi_E : E \to X, \quad \pi_F : F\to X, \quad \pi_{E \oplus F} :
    E \oplus F \to X, \quad t = \pi_{F^1} : E \oplus F \to E.
\end{equation*}
Let $x\in K^0_\GR(X).$ Then $\f^E(x)= \pi_E^*(x)\otimes
\La(\pi_E^*(E),s_E).$
Now we use that $t^*\pi_E^*(x) = \pi_{E \oplus F}^*(x)$ and
$t^*\La(\pi_E^*(E), s_E) = \La(\pi_{E \oplus F}^*(E), s_E \circ
t)$. Since $s_E \circ t + s_{F^1} = s_{E\oplus F},$ Equations
\eqref{eq.mult.Thom0} and \eqref{eq.mult.Thom} then give
\begin{equation*}
    \La(\pi_{E \oplus F}^*(E\oplus F), s_{E\oplus F})
    = \La(\pi_{E \oplus F}^*(E), s_E \circ t)
    \otimes \La (\pi_{E \oplus F}^*(F), s_{F^1}).
\end{equation*}
Putting together the above calculations we obtain
\begin{multline*}
    \f^{F^1} \f^E (x) = t^*(\f^E(x)) \otimes \La(t^*(F^1),
    s_{F^1}) \\ = t^*\pi_E^*(x) \otimes t^*(\Lambda (\pi_E^*(E),
    s_E)) \otimes \La( t^*(F^1), s_{F^1}) \\
    = \pi_{E \oplus F}^*(x) \otimes \La(\pi_{E \oplus F}^*E,
    s_E \circ t) \otimes \La (\pi_{E \oplus F}^*F, s_{F^1}) \\
    = \pi_{E \oplus F}^*(x) \otimes \La(\pi_{E \oplus F}^*(E
    \oplus F), s_{E\oplus F}) = \f^{E\oplus F}_1 (x).
\end{multline*}
The proof is now complete.
\end{proof}

We are now ready to formulate and prove the Thom isomorphism in
the setting of gauge-equivariant vector bundles.

\begin{teo}\label{isotomaobsh}\
Let $X \to B$ be a $\GR$-bundle and $F \to X$ a complex
$\GR$--equivariant vector bundle, then  $\f^{F} : K^{i}_{\GR}(X)
\to K^{i}_{\GR}(F)$ is an isomorphism.
\end{teo}

\begin{proof}\
Assume first that $F$ is a trivial bundle, that is, that $F = X
\times_B \cV$, where $\cV \to B$ is a complex, finite-dimensional
$\GR$--equivariant vector bundle. We continue to assume that $B$
is compact.

Let us denote by $\underline{\CC}:= B \times \CC$ the
$1$-dimensional $\GR$-bundle with the trivial action of $\GR$ on
$\CC$. Also, let us denote by $P(\cV\oplus \underline{\CC})$ the
projective space associated to $\cV \oplus \underline{\CC}$. As a
topological space, $P(\cV\oplus \underline{\CC})$ identifies with
the fiberwise one-point compactification of $\cV$. The embeddings
$\cV \subset P(\cV\oplus \underline{\CC})$ and $\cV \times_B  X
\subset  P(\cV\oplus \underline{\CC}) \times_B  X$ then gives rise
to the following natural morphism  (Equation \eqref{eq:conext})
\begin{equation*}
    j:K^0_\GR(\cV)\to K^0_\GR(P(\cV\oplus \underline{\CC})),
    \qquad j:K^0_\GR(  \cV \times_B  X)\to K^0_\GR(
    P(\cV\oplus \underline{\CC}) \times_B  X).
\end{equation*}

Let $X$ be compact and let $x \in K^0_\GR(P(\cV \oplus
\underline{\CC}) \times_B  X)$ be arbitrary. The fibers of the
projectivization $P(V\oplus \underline{\CC})$ are complex
manifolds, so we can consider the analytical index of the
correspondent family of Dolbeault operators over $P(V\oplus 1)$
with coefficients in $x$ (cf.~\cite[page~123]{At-Period}). This
index is an element of $K^0_\GR(X)$ by the results of \cite{NT1}.
Taking the composition with $j$
(cf.~\cite[page~122-123]{At-Period}) we get a family of mappings
$\a_{X}:K^0_\GR( \cV \times_B  X)\to K^0_\GR(X)$, having the
following properties:
\begin{enumerate}[(i)]
\item\ $\a_X$ is functorial with the respect to
$\GR$--equivariant morphisms;
\item\  $\a_X$ is a morphism of $K^0_\GR(X)$-modules;
\item\ $\a_{B}(\la_V)=1\in K^0_\GR(B)$.
\end{enumerate}

Let $X^+ := X \cup B$ be the fiberwise one-point compactification
of $X$. The commutative diagram
\begin{equation*}
        \xymatrix{
        0\to K^0_\GR(\cV \times_B  X)\ar[r]&K^0_\GR(\cV \times_B X^+)
        \ar[r]\ar[d]^{\a_{X^+}}   & K^0_\GR(\cV)\ar[d]^{\a_B}\\
        0\to K^0_\GR(X)\ar[r] &K^0_\GR(X^+)\ar[r]&K^0_\GR(B)}
\end{equation*}
with exact lines allows us to define $\a_X$ for $X$ non-compact.

Let $x\in K^0_\GR(X)$, then by (ii)
\begin{equation}\label{a6}
\a_X(\la_V x)=\a_X(\la_\cV) x=x,\qquad \a\f=\Id.
 \end{equation}

Let $q:= \pi_F : F = \cV\times_B X\to X$, $p: X\times_B \cV\to X$,
$q_1:\cV\times_B X\times_B \cV \to \cV$ (onto the first entry),
$q_2:X\times_B \cV \to  B$, and by $\wt y\in K^0_\GR(X\times_B
\cV)$ we denote the element, obtained from $y$ under the mapping
$X\times_B \cV \to \cV \times_B  X$, $(x,v)\mapsto (-v,x)$ (such
that $\cV\times_B X \times_B \cV\to \cV\times_B X \times_B  \cV$,
$(u,x,v)\mapsto (-v,x,u)$ is homotopic to the identity).

Let $y\in K^0_\GR(\cV \times_B X)$, then once again  by (i), (ii)
and then by (iii)
\begin{multline}\label{a7}
    \f(\a_X(y))= \pi_E^*\a_X(y)\otimes q^*\la_\cV
    =\a_{X\times_B \cV }(p^*_1 y)\otimes q^* \la_\cV=
    \a_{X\times_B \cV }(p^*_1y\otimes  q^*\la_V)\\
    =\a_{X\times_B \cV }(  y\boxtimes  \la_V)=
    \a_{X\times_B \cV }( \la_V  \boxtimes  \wt y)=
    \a_{X\times_B \cV }(q^*_1 \la_V  \otimes  \wt y)\\
    = \a_{X\times_B \cV }(q^*_1 \la_V)  \otimes  \wt y=
    q_2^* \a_{B}(\la_V) \otimes \wt  y =
    q_2^*(1) \otimes \wt  y  =\wt y\in K^0_\GR(X\times_B \cV),
\end{multline}
We obtain that $\f \circ \a_X$ is an isomorphism. Since $\a_X
\circ \f=\Id$, $\a_X$ is the two-sides inverse of $\f$ and the
automorphism $\f \circ \a_X$ is the identity.

The proof for a general (complex) $\GR$--equivariant vector bundle
$F\to X$ can be done as in \cite[p.~124]{At-Period}. However, we
found it more convenient to use the following argument. Embed
first $F$ into a trivial bundle $E = \cV\times_B X$. Let $\f^{1}$
and $\f^{2}$ be the Thom maps associated to the bundles $E \to F$
and $E \to X$. Then by Theorem \ref{thm:thom6} the diagram
\begin{equation*}
    \xymatrix{
    K^0_\GR(F)\ar[rr]^{\f^1}&&K^0_\GR(E)\\
    &K^0_\GR(X)\ar[ul]^{\f_F}\ar[ur]_{\f^2} & }
\end{equation*}
is commutative, while $\f^{2}$ is an isomorphism by the first part
of the proof. Therefore $\f_F$ is injective. The same argument
show that $\f^1$ is injective. But $\f^1$ must also be surjective,
because $\f^{2}$ is an isomorphism. Thus, $\f^1$ is an
isomorphism, and hence $\f^2$ is an isomorphism too.
\end{proof}

\section{Gysin maps}\label{sec:gysinmaps}

We now discuss a few constructions related to the Thom
isomorphism, which will be necessary for the definition of the
topological index. The most important one is the Gysin map. For
several of the constructions below, the setting of $\GR$-spaces
and even $\GR$-bundles
is
too general, and we shall have to consider {\em longitudinally
smooth} $\GR$--fiber bundles $\pi_X : X \to B$. The main reason
why we need longitudinally smooth bundles to define the Gysin map
is the same as in the definition of the Gysin map for embeddings
of smooth manifolds. We shall denote by $T_{\ve}X$ the {\em
vertical} tangent bundle to the fibers of $X \to B$. All tangent
bundles below will be {\em vertical} tangent bundles.

Let $X$ and $Y$ be \LS $\GR$-fiber bundles, $i: X\to Y$ be an
equivariant fiberwise embedding, and $p_T : T_{\ve}X \to X$ be the
vertical tangent bundle to $X.$ Assume $Y$ is equipped with a
$\GR$-invariant Riemannian metric and let $p_N : N_{\ve} \to X$ be
the fiberwise normal bundle to the image of $i$

Let us choose a function $\e: X \to (0,\infty)$ such that the map
of $N_{\ve}$ to itself \begin{equation*}
    n\mapsto \e\,\frac n{1+|n|}
\end{equation*}
is $\GR$-equivariant and defines a $\GR$-diffeomorphism $\F:
N_{\ve}\to W$ onto a bundle of open tubular neighborhoods
$W\supset X$ in $Y.$

Let $(N \oplus N)_{\ve} :=N_{\ve} \oplus N_{\ve}$. The embedding
$i : X\to Y$ can be written as a composition of two fiberwise
embeddings $i_1: X \to W$ and $i_2: W \to Y$.  Passing to
differentials we obtain
\begin{equation*}
    T_{\ve}X \ola{di_1}T_{\ve}W \ola{di_2}T_{\ve}Y\quad
    \text{and} \quad d\F:T_{\ve}N \to T_{\ve}W,
\end{equation*}
where we use the simplified notation $T_\ve N=T_\ve N_\ve$.

\begin{lem}
{\em (cf. \cite[page 112]{Friedrich})\/} The manifold $T_{\ve}N$
can be identified with $p^*_T (N\oplus N)_{\ve}$ with the help of
a $\GR$-equivariant dif\-fe\-o\-mor\-phism $\psi$ such that makes
the following diagram commutative
 \begin{equation*}
    \xymatrix{
    p^*_T(N\oplus N)_{\ve}\ar[d] && {T_\ve N} \ar[ll]^\psi\ar[d]\\
    T_{\ve}X\ar[dr]^{p_T} &  & N_{\ve}\ar[dl]_{p_N}\\
    &X.&}
 \end{equation*}
\end{lem}

\begin{proof}\
The vertical tangent  bundle $T_\ve N \to N_\ve$ and the vector
bundle
$$
p^*_N(T_{\ve}X)\oplus p^*_N(N_{\ve}) \to N_\ve
$$
are
isomorphic as $\GR$--equivariant vector bundles over $N_\ve$.

Indeed, a point of the total space $T_\ve N$ is a pair of the
form $(n_1, t + n_2),$ where both vectors are from the fiber over
the point $x\in X.$ Similarly, we represent elements
$p^*_T(N\oplus N)_{\ve}$ as pairs of the form $(t,n_1 + n_2).$ Let
us define $\psi$ by the equality $\psi (n_1, t+n_2)= (t,
n_1+n_2).$
\end{proof}

With the help of the relation $i\cdot (n_1,n_2)=(-n_2,n_1)$,
we can equip
\begin{equation*}
    p^*_T(N\oplus N)_{\ve}=p^*_T(N_{\ve})\oplus p^*_T(N_{\ve})
\end{equation*}
with a structure of a complex\label{ccc} manifold.
Then we can consider the Thom homomorphism
\begin{equation*}
    \f:K^0_\GR(T_\ve X)\to K^0_\GR(p^*_T(N\oplus N)_\ve).
\end{equation*}
Since $T_{\ve}W$ is an open $\GR$-stable subset of $T_{\ve}Y$ and
$di_2: T_{\ve}W\to T_{\ve}Y$ is a fiberwise embedding, by Equation
\ref{eq:conext}, there is the homomorphism
$(di_2)_*:K^0_\GR(T_\ve W)\to K^0_\GR(T_{\ve}Y). $

\begin{dfn}\label{oprgys}
Let $i: X\to Y$ be an equivariant embedding of $\GR$-bundles. The {\em
Gysin homomorphism\/} \label{ind:Gysinhomomorphism} is the mapping
\begin{equation*}
    i_!:K^0_\GR(T_{\ve}X)\to K^0_\GR(T_{\ve}Y),
    \qquad i_! = (di_2)_* \circ (d\F^{-1})^* \circ \psi^*
    \circ \f.
\end{equation*}
In other words, it is obtained by passage to $K$-groups in the
upper part of the diagram
\begin{equation*}
    \xymatrix{
    p^*_T(N\oplus N)_{\ve}\ar[d]_{q_T} &&T_\ve N
    \ar[ll]_\psi \ar[r]^{d\F}\ar[d]
    & T_{\ve}W \ar[r]^{di_2}\ar[dd]& T_{\ve}Y\ar[dd] \\
    T_{\ve}X \ar[dr]^{p_T} & & N_{\ve}\ar[dl]_{p_N} \ar[dr]^{\F} &
    & \\ & X\ar[rr]^{i_1} & & W \ar[r]^{i_2} & Y.}
\end{equation*}
 \end{dfn}

Another choice of metric and neighborhood $W$ induces the
homotopic map and (by the item 3 of Theorem~\ref{th:prophys}
below) the same homomorphism.

\begin{teo}[Properties of Gysin homomorphism]\label{th:prophys}
Let $i: X\to Y$ be a $\GR$-em\-bed\-ding.
\begin{enumerate}[\rm(i)]
\item\ $i_!$ is a homomorphism of $K^0_\GR(B)$-modules.
\item\ Let $i: X\to Y$ and $j: Y\to Z$ be two fiberwise
$\GR$-embeddings, then  $(j\circ i)_!=j_!\circ i_!.$
\item\ Let fiberwise embeddings $i_1: X\to Y$ and $i_2: X\to Y$ be
$\GR$-homotopic in the class of embeddings. Then
 $(i_1)_!=(i_2)_!.$
\item\ Let $i_!: X\to Y$ be a fiberwise $\GR$-diffeomorphism, then
 $
i_!=(di^{-1})^*.
 $
\item\
A fiberwise
embedding $i: X\to Y$ can be represented as a
compositions of embeddings $X$ in $N_{\ve}$
(as the zero section $s_0: x\to N$)
and $N_{\ve}\to Y$ by $i_2\circ \F: N_{\ve}\to Y.$ Then
 $
i_!=(i_2\circ \F)_! (s_0)_!.
 $
\item\ Consider the complex bundle $p^*_T (N_{\ve}\o\C)$ over
$T_{\ve}X$. Let us form the complex $\La(p^*_T(N_{\ve}\o\C),0):$
\begin{equation*}
    0\to\La^0(p^*_T(N_{\ve}\o\C))\ola{0}\dots
    \ola{0} \La^k(p^*_T(N_{\ve}\o\C))\to 0
 \end{equation*}
with noncompact support. If $a\in K^0_\GR(T_{\ve}X),$ then the
complex $$a\o \La(p^*_T(N_{\ve}\o\C),0)$$ has compact support and
defines an element of $K^0_\GR (T_{\ve}X)$. Then
\begin{equation*}
    (d i)^* i_! (a)=a\cdot \La(p^*_T(N_{\ve}\o\C),0),
\end{equation*}
where $di$ is the differential of the embedding $i$.
\item\
$i_!(x(di)^*y)=i_!(x)\cdot y,$ where $x\in K^0_\GR(T_{\ve}X)$ and
$y\in K^0_\GR(T_{\ve}Y)$.
\end{enumerate}
\end{teo}

\begin{proof}
(i)\ This follows from the definition of $i_!.$

(ii)\ To simplify the argument, let us identify the tubular
neighborhood with the normal bundle. Then $(j\circ  i)_!$ is the
composition
\begin{equation*}
    K^0_\GR(T_{\ve}X)\ola{\f} K^0_\GR(T_\ve N\oplus
    T_\ve N''_{\ve})\to K^0_\GR(T_{\ve}Z),
\end{equation*}
where $N'_{\ve}$ is the fiberwise normal bundle of $Y$ in $Z,$
$N''_{\ve}=N'_{\ve}|_X,$ and for the sum of tangent bundles to the
vertical normal bundles
 $T_\ve N\oplus
T_\ve N''_{\ve}$ is considered in the same way as on
page~\pageref{ccc}, that is, as a complex bundle over $T_{\ve}X$.
On the other hand, $j_!\circ i_!$ represents the composition
\begin{equation*}
    K^0_\GR(T_{\ve}X)\ola{\f} K^0_\GR(T_\ve N) \to
    K^0_\GR(T_{\ve}Y) \ola{\f} K^0_\GR(T_\ve N'_{\ve})
    \to K^0_\GR(T_{\ve}Z).
\end{equation*}
By the properties of $\f$, the following diagram is commutative
\begin{equation*}
    \xymatrix{
    K^0_\GR(T_\ve X)\ar[r]^{\f}\ar[dr]^\f & K^0_\GR(T_\ve N)\ar[r] \ar[d]^\f
    &K^0_\GR(T_\ve Y)\ar[d]^\f\\
    & K^0_\GR(T_\ve N\oplus T_\ve N'_\ve) \ar[r]^{\qquad\f} \ar[d]&
    K^0_\GR(T_\ve N'_\ve)\ar[d] \\
    & K^0_\GR(T_\ve Z) \ar@{=}[r]& K^0_\GR(T_\ve Z).}
\end{equation*}
This completes the proof of (ii).

(iii)\ The morphism $q_T$ depends only on the homotopy class of
the embeddings used to define it. The assertion thus follows from
the homotopy invariance of $K$-theory.

(iv)\ In this case $N=X,\,W=Y,\,\F=i,\,i_2=\Id_Y,$ and the formula
is obvious.

(v)\ This follows from (ii).

(vi)\ By definition,
 \begin{equation*}
(di)^* \circ i_! = (di_1)^* \circ (di_2)^* \circ (di_2)_*
\circ (d\f^{-1})^* \circ \psi^* \circ \f^* = (\psi \circ
d\F^{-1}\circ di_1)^*\circ \f,
 \end{equation*}
where $i_1: X\to W,\, i_2: W\to Y.$ Let
$(n_1,t+n_2)\in T_\ve N=p^*_N(T_\ve X)\oplus p^*_N(N_\ve)$,
where $n_1$ is the shift under the exponential mapping and
$t+n_2$ is a vertical tangent vector to $W.$
 If $d\F(n_1,t+n_2)$ is in $T_\ve X$, then  $n_1=n_2=0.$ Hence,
 \begin{equation*}
d\F^{-1}di_1(t)=(0,t+0),\quad
\psi\circ  d\F^{-1}\circ  di_1 (t)=(t,0+0).
 \end{equation*}
Therefore,
 $
\psi\circ  d\F^{-1}\circ  di_1 : T_\ve X \to p^*_T(N_\ve \oplus N_\ve)
 $
is the embedding of the zero section. Since
 $
\f(a)=a\cdot \La(q^*_T p^*_T (N_\ve \o \C), s_{p^*_T(N_\ve\o\C)}),
 $
it follows that
 $
(di)^* \circ i_!(a) = a \cdot \La(p^*_T(N_\ve\o\C),0).
 $

(vii)\ The mapping $di_1\circ  q_T \circ \psi\circ
d\F^{-1}:T_\ve W\to T_\ve W$ is homotopic to the identical mapping.
Hence,
\begin{equation*}
\begin{array}{l}
 i_!(x\cdot (di)^*y)=(di_2)_*(d\F^{-1})^*\psi^*\f(x\cdot (di)^*y)=\\
 \:=(di_2)_*(d\F^{-1})^*\psi^*[(q^*_T(x)\la_{p^*_T(N_\ve\o\C)})
(q^*_T(di)^*y)]=\\
 \:=(di_2)_*[(d\F^{-1})^*\psi^*(q^*_T(x)\la_{p^*_T(N_\ve\o\C)})
\underbrace{(d\F^{-1})^*\psi^*q^*_T(di_1)^*}_{\Id}(di_2)^*y]=\\
 \!=[(di_2)_*(d\F^{-1})^*\psi^*(q^*_T(x)\la_{p^*_T(N_\ve\o\C)})]\,
[(di_2)_*(di_2)^*y]=i_!(x)\cdot y.
         \end{array}
\end{equation*}
The proof is now complete.
\end{proof}

We shall need also the following properties of the Gysin map. If
$X = B$, the trivial \LS $\GR$-bundle, we shall identify $T_\ve X =
B$ and $T_\ve\cV = \cV\o\C$ for a real bundle $\cV\to B$.

\begin{teo}\label{th:prophys2}
Suppose that $\cV\to B$ is a $\GR$-equivariant {\em real} vector
bundle and that $X = B$. Then the mapping
\begin{equation*}
    i_!: K^0_\GR(B)=K^0_\GR(TX_\ve)\to
    K^0_\GR(T_{\ve}\cV)=K^0_\GR(\cV\o\C)
\end{equation*}
coincides with the Thom homomorphism $\f^{\cV\o\C}$.
\end{teo}

\begin{proof}\
The assertion follows from the definition of $i_!$. More
precisely, let $X=B\hookrightarrow \cV,\, N=\cV$ be the zero
section embedding. In the definition of the Thom isomorphism, $W$
can be chosen to be equal to the bundle $D_1$ of interiors of the
balls of radius 1 in $\cV$ with respect to an invariant metric. In
this case, the diagram from the definition of the Gysin
homomorphism \ref{oprgys} takes the following form
\begin{equation*}
\xymatrix{ \cV\o\C\ar[d]\ar[rr]^\Psi &&\cV\o\C\ar[r]^{d\F}\ar[d] &
D_1\o\C \ar[r]^{di_2}\ar[dd]& \cV\o\C\ar[dd] \\ T_\ve X =B \ar[dr]
& & \cV\ar[dl] \ar[dr]^{\F} & & \\
 & X=B\ar[rr]^{i_1} & & D_1 \ar[r]^{i_2} & \cV.
}
\end{equation*}
In our case $\Psi=\Id$ and $di_2\circ  d\F$ is homotopic to $\Id,$
since this map has the form
 $
v\o z \mapsto  (v\o z)/({1+|v\o z|}).
 $
Hence, $i_!=\f.$
\end{proof}

\begin{teo}\label{th:prophys3}
Suppose that $\cV'$ and $\cV''$ are $\GR$-equivariant $\R$-vector
bundles over $B$ and that $i: X\to \cV'$ is an embedding. Let $k :
X\to \cV'\oplus \cV''$, $k(x) =i(x) + 0.$ Let $\f$ be the Thom
homomorphism of the complex bundle
\begin{equation*}
T_\ve (\cV'\oplus \cV'')=\cV'\o\C \oplus \cV''\o\C
\ola{}T_\ve \cV'=\cV'\o\C.
\end{equation*}
Then the following diagram is commutative
\begin{equation*}
\xymatrix{ & K^0_\GR(T_\ve \cV')\ar[dd]^\f \\
K^0_\GR(T_\ve X)\ar[ur]^{i_!}\ar[dr]_{k_!}& \\ &
K^0_\GR(T_\ve (\cV'\oplus \cV'')). }
\end{equation*}
\end{teo}

\begin{proof}
To prove the statement, let us consider the embedding $i: X\to
\cV'$ and denote, as before, by $N_\ve$ the fiberwise normal
bundle and by $W$ the tubular neighborhood appearing in the
definition of the Thom isomorphism associated to $i$.  Then $N_\ve
\oplus \cV''$ is a fiberwise normal bundle for the embedding $k$
with the tubular neighborhood $W\oplus D (\cV''),$ where
$D(\cV'')$ is a ball bundle. If $a\in K^0_\GR (T_\ve X),$ then
\begin{equation*}
    k_!(a)=(di_2\oplus 1)_* \circ (d\F^{-1}\oplus 1)^*
    \circ (\psi\oplus 1)^* \circ \f^{N\oplus N \oplus
    \cV'' \oplus \cV''}(a).
\end{equation*}
We have $\f^{N\oplus N \oplus \cV'' \oplus \cV''}= \f^{N\oplus N}
\circ \f^{\cV'' \oplus \cV''},$ by Theorem \ref{thm:thom6}. Since
$a = a\cdot \underline {\C},$ where $\underline{\C}$ is the
trivial line bundle, we obtain
\begin{multline}
    k_!(a)=(di_2)_* \circ (d\F^{-1})^* \circ \Psi^* \circ
    \f^{N_\ve\oplus N_\ve} (a) \circ \f^{\cV'' \oplus \cV''}
    (\underline{\C}) \\ = i_!(a) \cdot
    \la_{T(\cV'\oplus \cV'')_\ve} = \f (i_!(a)).
\end{multline}
The proof is now complete.
\end{proof}

 \section{The topological index}\label{sec:topindex}

We begin with a ``fibered Mostow-Palais theorem''
that will be useful in defining the index.

\begin{teo}\label{teo:fibMosPalais}
Let $\pi_{X} : X \to B$ be a compact $\GR$-fiber bundle. Then
there exists a real $\GR$--equivariant vector bundle $\cV \to B$
and a fiberwise smooth $\GR$-embedding $X \to \cV$. After
averaging one can assume that the action of $\GR$ on $\cV$ is
orthogonal.
\end{teo}

\begin{proof}
Fix $b \in B$ and let $U$ be an equivariant trivialization neighborhood of $b$
for both $X$ and $\GR$. By the Mostow-Palais theorem,
there exists a representation of $\GR_b$ on a finite dimensional vector space
$V_b$ and a smooth $\GR_b$-equivariant embedding $i_b: X_b \to V_b$. This
defines an embedding
\begin{equation}\label{eq.embedding}
    \psi : \pi_X^{-1}(U_b) \simeq U_b \times X_b \to U_b \times V_b
\end{equation}
which is $\GR$-equivariant in an obvious sense.

We can cover $B$ with finitely many open sets $U_{b_j}$, as above,
corresponding to the points $b_j$, $j = 1, \ldots, N$. Denote by
$V_j$ the corresponding representations and by $\psi_j$ the
corresponding embeddings, as in Equation \eqref{eq.embedding}. Let
$W := \oplus V_j$. Also, Let $\phi_j$ be a partition of unity
subordinated to the covering by $U_j=U_{b_j}$. We define then
\begin{equation*}
    \Psi : = \oplus (\phi_j \circ \pi_X) \psi_j :
    X \to B \times W,
\end{equation*}
which is a $\GR$-equivariant embedding of $X$ into the trivial
$\GR$--equivariant vector bundle $\cV := B \times W$, as desired.
\end{proof}

Let us now turn to the definition of the topological index.
Let $X\to B$ be a compact, \LS $\GR$-bundle. {}From Theorem
\ref{teo:fibMosPalais} it follows that there exists  an
$\GR$--equivariant {\em real} vector bundle $\cV\to B$ and a
fiberwise smooth $\GR$-equivariant embedding $i:X\to \cV$. We can
assume that $\cV$ is endowed with an orthogonal metric and that
$\GR$ preserves this metric. Thus, the Gysin homomorphism
\begin{equation*}
    i_! : K^0_\GR(T_\ve X)\to K^0_\GR(T_\ve \cV)
    = K^0_\GR(\cV\o\C)
\end{equation*}
is defined  (see Section~\ref{subs:Thom}). Since
$T_\ve \cV=\cV\o\C$ is a complex vector bundle, we have the
following Thom isomorphism (see Section~\ref{subs:Thom}):
\begin{equation*}
    \f:K^0_\GR(B) \stackrel{\sim}{\longrightarrow} K^0_\GR(T_\ve \cV).
\end{equation*}

\begin{dfn}\label{dfi:topind}\ {\rm The {\em topological index\/}
is by definition the morphism:
 \begin{equation*}
    \ti^X_\GR:K^0_\GR(T_{\ve}X)\to K^0_\GR(B),\qquad
    \ti^X_\GR:=\f^{-1}\circ i_!.
 \end{equation*}
}
\end{dfn}

The topological index satisfies the following properties.

\begin{teo}\label{th:ax2} 
Let $X \to B$ be a \LS bundle and
$$\ti^X_\GR:K^0_\GR(T_{\ve}X)\to K^0_\GR(B)$$
be its associated topological index. Then
\begin{enumerate}[\rm(i)]
\item\ $\ti^X_\GR$ does not depend on the choice
of the $\GR$--equivariant vector bundle $\cV$ and on the embedding
$i: X\to \cV$.
\item\ $\ti^X_\GR$ is a $K^0_\GR(B)$-homomorphism.
\item\ If $X=B$, then the map
\begin{equation*}
    \ti^X_\GR : K^0_\GR(B)=K^0_\GR(T_{\ve} X)\to K^0_\GR(B)
\end{equation*}
coincides with $\Id_{K^0_\GR(B)}$.
\item\ Suppose $X$ and $Y$ are compact \LS $\GR$-bundles, $i: X
\to Y$ is a fiberwise $\GR$-embedding. Then the diagram
 \begin{equation*}
    \xymatrix{
    K^0_\GR(T_\ve X)\ar[rr]^{i_!}\ar[dr]_{\ti^X_\GR\:\:}
    && K^0_\GR(T_{\ve} Y)
    \ar[dl]^{\ti^Y_\GR}\\
    &K^0_\GR(B).&
}
 \end{equation*}
commutes.
 \end{enumerate}
\end{teo}

 \begin{proof}
To prove (i), let us consider two embeddings
 \begin{equation*}
i_1:X\to \cV',\qquad i_2:X\to \cV''
 \end{equation*}
into $\GR$--equivariant vector bundles. Denote by $j = i_1 + i_2$
the induced embedding $j: X\to \cV'\oplus \cV''$. It is
sufficient to show that $i_1$ and $j$ define the same topological
index. Let us define a homotopy of $\GR$-embeddings by the
formula
 \begin{equation*}
    j_s(x)=i_1(x)+s\cdot i_2(x)\::\:X\to \cV'\oplus \cV'',
    \quad 0\le s\le 1.
 \end{equation*}
Then, by Theorems~\ref{th:prophys}(iii) and \ref{th:prophys3}, the
indices for $j$ and $j_0$ coincide. Let us show now that $j_0=i_1
+ 0$ and $i_1$ define the same topological indexes. For this
purpose consider the diagram
 \begin{equation*}
\xymatrix{ &K^0_\GR(T_{\ve}X)\ar[dl]_{(i_1)_!}\ar[dr]^{(j_0)_!}&\\
K^0_\GR(T_\ve \cV')\ar[rr]^{\f_2} && K^0_\GR(T_\ve (\cV'\oplus \cV'')), \\
&K^0_\GR(B)\ar[ul]^{\f_1}\ar[ur]_{\f_3}& }
 \end{equation*}
where $\f_i$ are the corresponding Thom homomorphisms. The upper
triangle is commutative by Theorem \ref{th:prophys2}.2, and the
lower is commutative by Theorem \ref{thm:thom6}. Hence
$\f_1^{-1}\circ (i_1)_!=\f_3^{-1}\circ  (j_0)_1$ as desired.

(ii) follows from \ref{prop:thom3} and \ref{th:prophys}(i).

Property (iii) follows from the definition of the index and from
\ref{th:prophys2}.

To prove (iv), let us consider the diagram
 \begin{equation*}
\xymatrix{
X\ar[rr]^i\ar[dr]_{j\circ i}&&Y\ar[dl]^j\\
&\cV.&
}
 \end{equation*}
We now use \ref{th:prophys}(ii). This gives the commutative
diagram
 \begin{equation*}
    \xymatrix{ K^0_\GR(T_{\ve}X)\ar[rr]^{i_!}\ar[dr]_{(j\circ i)_!\:
    }&&K^0_\GR(T_{\ve}Y)\ar[dl]^{j_!}\\ &K^0_\GR(T_\ve \cV)& }
 \end{equation*}
or
 \begin{equation*}
\xymatrix{ K^0_\GR(T_{\ve}X)\ar[rr]^{i_!}\ar[dr]^{(j\circ
i)_!\:}\ar[ddr]_{\ti^X_\GR\:}
&&K^0_\GR(T_{\ve}Y)\ar[dl]_{j_!}\ar[ddl]^{\,\ti^Y_\GR}\\
&K^0_\GR(T_\ve \cV)&\\ &K^0_\GR(B).\ar[u]_\f& }
 \end{equation*}
This completes the proof.
 \end{proof}

We now investigate the behavior of the topological index
with respect to fiber products of bundles of compact groups.

 \begin{teo}\label{th:ax7}
Let $\pi: P\to X$ be a principal right $\H$-bundle with a left
action of $\GR $ commuting with $\H$. Suppose $F$ is a \LS $(\GR
\times \H)$-bundle. Let us denote by $Y$ the space $P\times_\H F$.
Let $j: X'\to X$ and $k: F'\to F$ be fiberwise $\GR $- and $ (\GR
\times \H)$-embeddings, respectively. Let $\pi': P'\to X'$ be the
principal $\H$-bundle induced by $j$ on $X'$. Assume that $Y' :=
P'\times_\H F'$. The embeddings $j$ and $k$ induce $\GR
$-embedding $j*k: Y'\to Y$. Then the diagram
\begin{equation*}
\xymatrix{ K^0_\GR(T_{\ve}X)\otimes_{K^0_\GR (B)} K^0_{\GR \times
\H}(T_{\ve}F) \ar[r]^{\qquad\qquad \g} & K^0_\GR(T_{\ve}Y)\\
K^0_\GR(T_{\ve}X')\otimes_{K^0_\GR (B)} K^0_{\GR \times
\H}(T_\ve F')\ar[r]^{\qquad \qquad \g } \ar[u]_{j_!\o k_!} &
K^0_\GR(T_\ve Y')\ar[u]_{(j*k)_!} }
 \end{equation*}
is commutative.
 \end{teo}

Let us remark that there in the statement of this theorem
there is no compactness assumption on $X, X', F,$ and $F'$, since
there is no compactness assumption in the definition of the Gysin
homomorphism. This is unlike in the definition of the topological
index where we start with a compact $\GR$-bundle $X \to B$.

 \begin{proof}
{Let us use the definition of $\g$:}
 \begin{equation*}
\xymatrix{ K^0_\GR(T_{\ve}X)\otimes K^0_{\GR \times
\H}(T_{\ve}F)\ar[r]\ar@{}[dr]|{\fbox{1}}& K^0_\GR(T_{\ve}X)\o
K^0_{\GR \times \H}(P\times T_{\ve}F)\ar[r] ^{\qquad \qquad \qquad
\cong} \ar@{}[dr]|{\fbox{2}} &\\ K^0_\GR(T_{\ve}X')\otimes
K^0_{\GR \times \H}(T_\ve F')\ar[r]\ar[u]^{j_!\o k_!} &
K^0_\GR(T_{\ve}X')\o K^0_{\GR \times \H}(P'\times
T_\ve F')\ar[u]^{\e}\ar[r] ^{\qquad \qquad \qquad \cong} & }
 \end{equation*}
\begin{equation*}
\xymatrix{
\cong 
K^0_\GR(T_{\ve}X)\o K^0_\GR(P\times_\H T_{\ve}F)\to
\ar@{}[dr]|{\fbox{3}} & \hspace{-3em} K^0_\GR(\pi^*_1 T_{\ve}X)\o
K^0_\GR(P\times_\H T_{\ve}F)\to
\\
\cong  K^0_\GR(TX')\o K^0_\GR(P'\times_\H T_\ve F')\to
\ar[u]^{\b} & \hspace{-3em} K^0_\GR((\pi'_1)^* T_{\ve}X')\o
K^0_\GR(P'\times_\H T_\ve F')\to\ar[u]^{\a} }
\end{equation*}
\begin{equation}\label{axio1}
 \begin{array}{c}
\to K^0_\GR((\pi^*_1  T_{\ve}X)\times(P\times_\H T_{\ve}F)) \to \\
         \fbox{4} \\
\to K^0_\GR((\pi'_1)^*  T_{\ve}X'\times(P'\times_\H T_\ve F'))
\to
 \end{array}
\end{equation}
\smallskip
 \begin{equation*}
 \begin{array}{c}
    \to K^0_\GR(\pi^*_1 T_{\ve}X\oplus (P\times_\H
    T_{\ve}F))=K^0_\GR(T_{\ve}Y)\\ \qquad \qquad \fbox{4}\qquad
    \qquad\qquad \qquad \qquad \qquad \uparrow\, (j*k)_! \\ \to
    K^0_\GR((\pi'_1)^* T_{\ve}X'\oplus (P'\times_\H
    T_\ve F'))=K^0_\GR(TY'_{\ve}),
 \end{array}
 \end{equation*}
where projections $\pi_1 : Y = P\times_\H F\to X$ and
$\pi'_1:Y'=P\times_\H F'\to X'$ are defined as above. Here we use
the isomorphism $K^0_{\GR\times\H}(P\times W)\cong
K^0_\GR(P\times_{\H}W)$ for a free $\H$-bundle $P$ (see
Theorem~\ref{lem:superinduction}). Let us remind the diagram,
which was used for the definition of the Gysin homomorphism of an
embedding $j: X'\to X$:
 \begin{equation*}
\xymatrix{ p^*_T(N_{X'}\oplus N_{X'})_\ve \ar[d]^{q^{X'}_T}
&&T_\ve N_{X'}\ar[ll]_\psi\ar[r]^{d\F_{X'}} \ar[d] &
T_\ve W_{X'} \ar[r]^{dj_2}\ar[dd]& T_{\ve}X \ar[dd] \\ T_{\ve}X'
\ar[dr]^{p_T} & & N_{X',{\ve}}\ar[dl]_{p_{N_{X',{\ve}}}\:}
\ar[dr]^{\F_{X'}}&& \\ & X'\ar[rr]^{j_1} & & W_{X'} \ar[r]^{j_2} &
X.  }
\end{equation*}
From the similar diagrams for $k_!$ and $ (j * k) _!$ and the
explicit form of these maps, it follows that the square \fbox{4}
in (\ref {axio1}) is commutative if, and only if, $\a$ has the
following form:
 \begin{equation*}
 \begin{array}{l}
\a(\s\o \rho)=(\pi^*_1)\Bigl\{
(dj_2)_*\,(d\F^{-1}_{X'})^*\,\psi^*_{X'}\Bigr\}\circ\f^S(\s)\o\\
\qquad \o(\pi^*j_2\times_\H  dk_2)_*\,\Bigl((\pi^*\F_{X'}\times_\H
d\F_{F'})^{-1}\Bigr)^*\,(1\times_\H\psi_{F'})^*\,\f^R(\rho),
 \end{array}
 \end{equation*}
where $S$ and $T$ are bundles of the form
 \begin{equation*}
 \begin{array}{cccc}
 &\pi^*_1\,\Bigl((p^{X'}_T)^*\{N_{X'}\oplus N_{X'}\}\Bigr)&&
\pi^* N_{X'}\times_\H (p^{F'}_T)^*\,(N_{F'}\oplus N_{F'})\\
S:&\qquad \downarrow\,(\pi'_1)^*q^{X'}_T & R:&\qquad
\downarrow\,(\pi')^*(p_{N_{X'}})\times_\H q^{F'}_T \\
&(\pi'_1)^*\,(T_{\ve}X'), && \pi^*\,X'\times_\H
T_\ve F'=P'\times_\H T_\ve F'.
 \end{array}
 \end{equation*}
Hence the square \fbox{3} in (\ref {axio1}) is commutative iff the
homomorphism $\b$ has the form
 \begin{equation*}
 \begin{array}{l}
\b(\t\o\rho)=j_!(\t)\o \\
\:\o(\pi^*j_2\times_\H dk_2)_*
\Bigl((\pi^*\F_{X'}\times_\H
d\F_{F'})^{-1}\Bigr)^*\,(1\times_\H\psi_{F'})^*\,\f^R(\rho),
 \end{array}
  \end{equation*}
 \begin{equation*}
\t\in K^0_\GR(TX'),\qquad \rho\in K^0_\GR(P'\times_\H TF').
 \end{equation*}
In turn, the square \fbox{2} in (\ref {axio1})
is commutative iff the homomorphism $\e$ has the form
 \begin{equation*}
 \begin{array}{l}
\e(\t\o\delta)=j_!(\t)\o \\
\:\o(\pi^*j_2\times_\H dk_2)_*
\Bigl((\pi^*\F_{X'}\times_\H
d\F_{F'})^{-1}\Bigr)^*\,(1\times_\H\psi_{F'})^*\,\f^{\widetilde
R}_\C(\delta),
\end{array}
  \end{equation*}
 \begin{equation*}
\t\in K^0_\GR(TX'),\qquad \delta\in K^0_{G\times \H}(P'\times TF'),
 \end{equation*}
where $\widetilde R$ is the following bundle:
 \begin{equation*}
 \begin{array}{cc}
&\pi^* N_{X'}\times (p^{F'}_T)^*\,(N_{F'}\oplus N_{F'})\\
\widetilde R:\quad&\qquad \downarrow\,(\pi')^*(p_{N})\times
q^{F'}_T \\
&P'\times TF'.\qquad \qquad
 \end{array}
 \end{equation*}
Suppose $\delta=[{\underline \C}]\widehat\o\omega$, where
$[{\underline \C}]\in K^0_{\GR \times \H}(P')$, ${\underline \C}$ is
the one-dimensional trivial bundle and $\omega\in K^0_{\GR \times \H}(TF')$.
Then
\begin{eqnarray*}
\e(\t\o\delta)&=&j_!(\t)\o\Bigl\{\pi^*(j_2)_*(\F^{-1}_{X'})^*
[{\underline \C}]\widehat\o k_!(\omega)\Bigr\}=\\
&=& j_!(\t)\o\Bigl\{[{\underline \C}]\widehat\o k_!(\omega)\Bigr\}.
\end{eqnarray*}
Since the map $K^0_{\GR \times \H}(TF)\to K^0_{\GR \times \H}(P\times TF)$ (as well
as the lower line in (\ref{axio1})) has the form
 $
\omega\mapsto [{\underline \C}]\widehat\o\omega ,
 $
we have proved the commutativity of \fbox {1} in (\ref{axio1}).
 \end{proof}

From this theorem we obtain the following corollary.

 \begin{cor}\label{cor:commaxio}
Let $M$ be a compact smooth $H$-manifold, let $\H=B\times H$, and
let $P$ be a principal \LS $\H$-bundle over $X$ carrying also an
action of $\GR$ commuting with the action of $\H$. Also, let $X\to
B$ be a compact \LS $\GR$-bundle. Let $Y:=P\times_H M\to X$ be
associated \LS $\GR$-bundle. Taking $F=B\times M$, we define $T_M
Y:=T_F Y$. Then $T_M Y$ is a $\GR$-invariant real subbundle of
$T_{\ve}Y$ and $T_M Y=P\times_H TM$. Let $j: X'\to X$ be a
fiberwise $\GR$--equivariant embedding and let $k: M'\to M$ be an
$H$-embedding. Denote by $\pi': P'\to X'$ the principal
$\H$-bundle induced by $j$ on $X'$ and assume that $Y':=P'\times_H
M'$. The embeddings $j$ and $k$ induce $\GR $-embedding $j*k:
Y'\to Y$. Then the diagram
\begin{equation*}
\xymatrix{ K^0_\GR(T_{\ve}X)\otimes K^0_{H}(TM) \ar[r]^{\qquad \g} &
K^0_\GR(T_{\ve}Y)\\ K^0_\GR(T_{\ve}X')\otimes K^0_H(TM')\ar[r]^{\qquad \g
} \ar[u]_{j_!\o k_!} & K^0_\GR(T_\ve Y')\ar[u]_{(j*k)_!} }
 \end{equation*}
is commutative.
\end{cor}

\providecommand{\bysame}{\leavevmode\hbox
to3em{\hrulefill}\thinspace}
\providecommand{\MR}{\relax\ifhmode\unskip\space\fi MR }
\providecommand{\MRhref}[2]{%
  \href{http://www.ams.org/mathscinet-getitem?mr=#1}{#2}
} \providecommand{\href}[2]{#2}

\end{document}